\documentstyle[amsfonts, amssymb, latexsym,12pt]{amsart}

\newtheorem{theorem}{Theorem}[section]
\newtheorem{lemma}[theorem]{Lemma}
\newtheorem{proposition}[theorem]{Proposition}
\newtheorem{corollary}[theorem]{Corollary}
\newtheorem{definition}[theorem]{Definition}

\def\Z{{\mbox{\rm\kern.25em
\vrule width.03em height0.57ex depth0ex
\kern.033em
\vrule width.03em height1.52ex depth-0.96ex \kern-.338em Z}}}
\def\R{{\mbox{\rm I\kern-.22em R}}}
\def\N{{\mbox{\rm I\kern-.22em N}}}
\def\P{{\bf P}}

\def\D{{\bf D}}
\def\T{{\bf T}}
\def\supp{{\rm supp}}
\def\size{{\rm size}}

\def\energy{{\rm energy}}

\def\D{{\cal{D}}}
\def\J{{\cal{J}}}

\def\dist{{\rm dist}}

\def\111{\gamma}

\def\be#1{\begin{equation}\label{#1}}
\def\bas{\begin{align*}}
\def\eas{\end{align*}}
\def\bi{\begin{itemize}}
\def\ei{\end{itemize}}
\newenvironment{proof}{\noindent {\bf Proof} }{\endprf\par}
\def \endprf{\hfill  {\vrule height6pt width6pt depth0pt}\medskip}
\def\emph#1{{\it #1}}

\title{Generalizations of the Carleson-Hunt theorem I.\\ The classical 
singularity case}

\author{Xiaochun Li}
\address{Department of Mathematics, University of Illinois at Urbana-Champaign,
Urbana, IL 61801}
\email{xcli@@math.uiuc.edu}

\address{Current Address: School of Mathematics, IAS, Princeton, NJ 08540}
\email{xcli@@math.ias.edu}

\author{Camil Muscalu}
\address{Department of Mathematics, Cornell University, Ithaca, NY 14853}
\email{camil@@math.cornell.edu}

\begin{document}

\begin{abstract}
In this article, we prove $L^p$ estimates for a general maximal
operator, which extend both the classical Coifman-Meyer \cite{meyerc} and 
Carleson-Hunt \cite{carleson}, \cite{hunt} theorems in harmonic analysis.
\end{abstract}

\maketitle

\section{Introduction}

This article is the first in a sequel of papers whose aim is to present 
several generalizations of the celebrated Carleson-Hunt theorem in Fourier
analysis.

The maximal Carleson operator is the sub-linear operator defined by

\begin{equation}\label{carleson-def}
Cf(x):= \sup_{N\in\R} 
\left|\int_{\xi<N}\widehat{f}(\xi) e^{2\pi i x\xi} d\xi \right|,
\end{equation}
where $f$ is a Schwartz function on $\R$ and the Fourier transform is defined 
by

\begin{equation}
\widehat{f}(\xi):= \int_{\R} f(x) e^{-2 \pi i x\xi} dx.
\end{equation}
The following result of Carleson and Hunt \cite{carleson}, \cite{hunt}
 is a classical 
theorem in Fourier analysis.

\begin{theorem}\label{carleson-thm}
The operator $C$ maps $L^p\rightarrow L^p$ boundedly, for every $1<p<\infty$.
\end{theorem}
This statement, in the particular weak type $L^2\rightarrow L^{2,\infty}$ 
special case, was the main ingredient in the proof of Carleson's fameous
theorem which states that the Fourier series of a function in $L^2(\R/\Z)$
converges pointwise almost everywhere.

For $n\geq 1$, let now consider $m(=m(\xi))$ in $L^{\infty}(\R^n)$ a bounded
function, smooth away from the origin and satisfying

\begin{equation}\label{mihlin}
|\partial^{\alpha}m(\xi)|
\lesssim\frac{1}{|\xi|^{|\alpha|}}
\end{equation}
for sufficiently many multi-indices $\alpha$. Denote by $T_m$ the $n$-linear 
operator defined by

\begin{equation}\label{tm}
T_m(f_1,...,f_n)(x):=\int_{\R^n}m(\xi) \widehat{f_1}(\xi_1)...
\widehat{f_n}(\xi_n) e^{2 \pi i x (\xi_1+...+\xi_n)} d\xi
\end{equation}
where $f_1,...,f_n$ are Schwartz functions on the real line $\R$. The 
following statement of Coifman and Meyer \cite{meyerc} is also a 
classical theorem in analysis.

\begin{theorem}\label{cm-thm}

$T_m$ maps $L^{p_1}\times...\times L^{p_n}\rightarrow L^p$ boundedly, as 
long as $1<p_1,...,p_n\leq\infty$, 
$\frac{1}{p_1}+...+\frac{1}{p_n}=\frac{1}{p}$ and $0<p<\infty$.
\end{theorem}

Now, for $N\in\R^n$ and $m$ as before satisfying (\ref{mihlin}), denote 
by $\tau_N m(\xi):= m(\xi-N)$ the translated symbol
and by $C_m$ the maximal operator defined by

\begin{equation}\label{carlesonm-def}
C_m(f_1,...,f_n)(x):=
\sup_{N\in\R^n}
\left|
\int_{\R^n}\tau_N m(\xi) \widehat{f_1}(\xi_1)...
\widehat{f_n}(\xi_n) e^{2 \pi i x (\xi_1+...+\xi_n)} d\xi
\right|
\end{equation}
where as before $f_1,...,f_n$ are Schwartz functions on $\R$.

The purpose of the present paper is to study the $L^p$ boundedness properties
of this Carleson type operator $C_m$. Our main theorem is the following.

\begin{theorem}\label{main}

$C_m$ maps $L^{p_1}\times...\times L^{p_n}\rightarrow L^p$ boundedly, as 
long as $1<p_1,...,p_n\leq\infty$, 
$\frac{1}{p_1}+...+\frac{1}{p_n}=\frac{1}{p}$ and $0<p<\infty$.
\end{theorem}

Clearly, Theorem \ref{main} contains both Coifman-Meyer theorem and 
Carleson-Hunt theorem as special cases.

To motivate the introduction of this operator, we should mention that a 
simplified variant of it appeared recently in connection to the so called
bi-Carleson operator studied in \cite{mtt:discrete-bicarleson} and
\cite{mtt:bicarleson}. This is the operator defined by the following formula

\begin{equation}\label{bicarleson-def}
T(f_1,f_2)(x):= \sup_{N\in\R}
\left|
\int_{\xi_1<\xi_2<N}\widehat{f_1}(\xi_1)\widehat{f_2}(\xi_2)
e^{2\pi i x (\xi_1 +\xi_2)} d\xi_1 d \xi_2
\right|
\end{equation}
and the following estimates are known about it \cite{mtt:discrete-bicarleson},
\cite{mtt:bicarleson}.

\begin{theorem}\label{bicarleson-thm}

$T$ maps $L^{p_1}\times L^{p_2}\rightarrow L^p$ boundedly, as long as
$1<p_1, p_2\leq\infty$, $\frac{1}{p_1}+\frac{1}{p_2} = \frac{1}{p}$ and
$\frac{2}{3}<p<\infty$.
\end{theorem}
If one denotes by $m(\xi_1,\xi_2):=\chi_{\{\xi_1<\xi_2<0\}}$, then one 
observes that 

\begin{equation}\label{simpler}
T(f_1,f_2)(x):= \sup_{N\in\R}
\left|
\int_{\R^2}\tau_{(N,N)}m(\xi_1, \xi_2)\widehat{f_1}(\xi_1)\widehat{f_2}(\xi_2)
e^{2\pi i x (\xi_1 +\xi_2)} d\xi_1 d \xi_2
\right|.
\end{equation}
In \cite{mtt:bicarleson}, the authors decomposed the symbol $m(\xi_1,\xi_2)$
as

$$m(\xi_1,\xi_2)= u(\xi_1,\xi_2) + v(\xi_1,\xi_2) + w(\xi_1,\xi_2)$$
where $u(\xi_1,\xi_2)$ was a symbol singular only along $\xi_2=0$, 
$v(\xi_1,\xi_2)$ was a symbol singular only along $\xi_1=\xi_2$, while 
$w(\xi_1,\xi_2)$
was a symbol singular only at the origin and satisfying (\ref{mihlin}).
Consequently, the operator $T$ could be estimated by a sum of three distinct
operators $U + V + W$, each corresponding to the symbols $u, v$ and $w$ 
respectively. In \cite{mtt:discrete-bicarleson}, \cite{mtt:bicarleson}
the operators $U$ and $V$ have been studied carefully, while the estimates
for $W$ followed from the main Theorem \ref{main} of this paper
(in fact, the operator $W$ is simpler than $C_m$ in 
(\ref{carlesonm-def}) since the translations in (\ref{simpler}) are made 
only along the line $\xi_1=\xi_2$ and not in the whole plane $\R^2$).

For the simplicity of our exposition and also for the reader's convenience
we chose to present the proof of our main Theorem \ref{main} 
in the particular case $n=2$. However, it will be clear from the proof
that its extension to the $n$-sub-linear case is straightforward.
While the current article is essentially selfcontained, we adopt the same 
strategy as in \cite{mtt:fourierbiest}, \cite{mtt:bicarleson} and will mark as
``standard'' any results that are well understood by now in this framework,
as in \cite{carleson}, \cite{fefferman}, \cite{grafakosli}, \cite{laceyli},
\cite{laceyt2}, \cite{mtt:fourierbiest}, \cite{mtt:bicarleson},
\cite{thiele}, etc.

The authors have been partially supported by NSF. The second author was also 
partially supported by an Alfred P. Sloan Research Fellowship. Both of the 
authors would like to thank Michael Lacey, Terry Tao and Christoph Thiele
for valuable conversations.

\section{Notation}

In this section we set out some general notations used throughout the paper.
We will write $A\lesssim B$ to denote the statement that $A\leq CB$ for some 
large constant $C>0$ and $A\sim B$ to denote the statement that 
$A\lesssim B\lesssim A$. Given an arbitrary interval $I$, we denote by $|I|$
the measure of $I$ and by $cI$ $(c>0)$ the interval with the same center
as $I$ but $c$ times its length. Given an interval $I$ we also denote
the approximate cutoff function $\tilde{\chi_I}$ by

$$\tilde{\chi_I}(x):= 
\left(
1+\frac{\dist(x, I)}{|I|}
\right)^{-10}.
$$
We say that a smooth function $\Phi_I$ is a bump adapted to $I$ if and only if
the following inequalities hold

\begin{equation}
|\Phi_I^l(x)|\leq
C_{l,m}\frac{1}{|I|^l}\frac{1}{\left(1+\frac{\dist(x,I)}{|I|}\right)^m}
\end{equation}
for every integer $m\in\N$ and for sufficiently many derivatives $l\in\N$.
If $\Phi_I$ is a bump adapted to $I$, we say that $|I|^{-1/2}\Phi_I$
is an $L^2$ normalized bump adapted to $I$.

\section{A ``translation invariant'' Littlewood-Paley decomposition 
of the symbol $m$}

In this section we describe a decomposition of the symbol $m$, which is well 
adapted to its future translations over vectors $N\in\R^2$.

Let $M>0$ be a big integer which will be fixed throughout the paper. For
$j=1,2,...,M$ consider Schwartz functions $\psi_j$ so that
$\supp \widehat{\psi_j}\subseteq \frac{10}{9}[j-1,j]$, $\psi_j=1$ on
$[j-1,j]$ and for $j=-M,-M+1,...,-1$ consider Schwartz functions $\psi_j$
so that $\supp \widehat{\psi_j}\subseteq \frac{10}{9}[j,j+1]$ and
$\psi_j=1$ on $[j,j+1]$.

For $\lambda>0$ and $\psi$ a Schwartz function, we denote by

$$D_{\lambda}^p\psi(x):= \lambda^{-1/p}\psi(\lambda^{-1}x)$$
the dilation operator which preserves the $L^p$ norm of $\psi$,
for $1\leq p\leq\infty$. Define the new symbol $\tilde{m}(\xi_1,\xi_2)$ by the 
formula

\begin{equation}\label{mtilda}
\tilde{m}(\xi_1,\xi_2):=
\sum_{\max(|j_1|, |j_2|)=M}
\int_{\R}
\left(D_{2^k}^{\infty}\widehat{\psi_{j_1}}(\xi_1)
D_{2^k}^{\infty}\widehat{\psi_{j_2}}(\xi_2)\right) dk.
\end{equation}
Clearly, by construction, $\tilde{m}$ is a bounded symbol, smooth away from
the origin and satisfying (\ref{mihlin}). Also, things can be arranged so that
$|\tilde{m}(\xi_1,\xi_2)|\geq c_0>0$ for every $(\xi_1,\xi_2)\in\R^2$ where
$c_0$ is a universal constant. This $\tilde{m}$ should be understood as being
essentially a decomposition of unity in frequency space, 
into a series of smooth
functions supported on rectangular annuli.

Then, 
we write $m(\xi_1,\xi_2)$ as

$$m(\xi_1,\xi_2)=
\frac{m(\xi_1,\xi_2)}{\tilde{m}(\xi_1,\xi_2)}\cdot
\tilde{m}(\xi_1,\xi_2):=\tilde{\tilde{m}}(\xi_1,\xi_2)\cdot
\tilde{m}(\xi_1,\xi_2)$$

\begin{equation}\label{tilde2}
=\sum_{\max(|j_1|, |j_2|)=M}
\int_{\R}
\left(\tilde{\tilde{m}}(\xi_1,\xi_2)
D_{2^k}^{\infty}\widehat{\psi_{j_1}}(\xi_1)
D_{2^k}^{\infty}\widehat{\psi_{j_2}}(\xi_2)\right) dk
\end{equation}
and observe that $\tilde{\tilde{m}}(\xi_1,\xi_2)$ has the same properties
as $m(\xi_1,\xi_2)$.

Fix now $j_1, j_2$ with $\max(|j_1|, |j_2|)=M$ and $k\in\R$. By taking 
advantage of the fact that $\tilde{\tilde{m}}$ satisfies (\ref{mihlin}),
one can write it on the support of 
$D_{2^k}^{\infty}\widehat{\psi_{j_1}}\otimes
D_{2^k}^{\infty}\widehat{\psi_{j_2}}$ as a double Fourier series and this 
allows us to decompose the corresponding inner term in (\ref{tilde2}) as

\begin{equation}\label{desc}
\sum_{n_1,n_2\in\Z}
C^{j_1,j_2}_{k,n_1,n_2}
\left(
D_{2^k}^{\infty}\widehat{\psi_{j_1}}(\xi_1)
e^{2\pi i n_1 \frac{9}{10}2^{-k}\xi_1}
\right)
\left(
D_{2^k}^{\infty}\widehat{\psi_{j_2}}(\xi_2)
e^{2\pi i n_2 \frac{9}{10}2^{-k}\xi_2}
\right)
\end{equation}
where

\begin{equation}\label{bound}
|C^{j_1,j_2}_{k,n_1,n_2}|\lesssim
\frac{1}{(1+|n_1|)^{1000}}
\frac{1}{(1+|n_2|)^{1000}}
\end{equation}
for every $n_1, n_2\in\Z$, uniformly for $k\in\R$.

The following Lemma will play a crucial role in our further decomposition

\begin{lemma}\label{lr}
Let $I=[a_I, b_I]$ be an interval so that $1\leq |I|\leq 2$ and $\Phi_I$ be a 
smooth function supported on it. Let also $\#\in (0,1)$ be fixed and much 
smaller than $1/M$ and define $\lambda_{\text{left}}$ and 
$\lambda_{\text{right}}$ by

\begin{equation}\label{l}
\lambda_{\text{left}}(x):=\frac{1}{\#}
\int_1^{1+\#}\chi_{[-\alpha, 0]}(x) d\alpha
\end{equation}
and

\begin{equation}\label{r}
\lambda_{\text{right}}(x):=\frac{1}{\#}
\int_1^{1+\#}\chi_{[0, \alpha]}(x) d\alpha.
\end{equation}
Then, there exist bump functions $\Phi_I^{\text{left}}$ and 
$\Phi_I^{\text{right}}$ adapted to I so that

\begin{equation}\label{*l}
\Phi_I^{\text{left}}\ast \lambda_{\text{left}} = \Phi_I
\end{equation}
and

\begin{equation}\label{*r}
\Phi_I^{\text{right}}\ast \lambda_{\text{right}} = \Phi_I.
\end{equation}
Moreover, they have the additional properties that 
$\supp \Phi_I^{\text{left}}\subseteq (-\infty, b_I]$ and
$\supp \Phi_I^{\text{right}}\subseteq [a_I, +\infty)$.
\end{lemma}

\begin{proof}
Since $\lambda_{\text{left}}$ is compactly supported, it follows that its 
Fourier transform $\widehat{\lambda_{\text{left}}}$ is smooth and
$|\widehat{\lambda_{\text{left}}}|\lesssim 1$. Moreover, we claim that
$\widehat{\lambda_{\text{left}}}(\xi)\neq 0$ for every $\xi\in\R$ and that
$|\widehat{\lambda_{\text{left}}}(\xi)|\sim \frac{1}{|\xi|}$ for
$|\xi|\geq 1$.

Clearly, $\widehat{\lambda_{\text{left}}}(\xi)\neq 0$ if $|\xi|$ is small
enough, say $|\xi|< \frac{1}{1000}$. On the other hand, for 
$|\xi|\geq \frac{1}{1000}$, one can write

\begin{equation}
\widehat{\lambda_{\text{left}}}(\xi)=
\frac{1}{\#}\int_1^{1+\#}
\left(\int_{\R}\chi_{[-\alpha, 0]}(x) e^{-2\pi i x\xi} dx\right) d\alpha
=\frac{1}{2\pi i \xi}\frac{1}{\#}\int_1^{1+\#}
(e^{2\pi i \alpha\xi}-1) d\alpha.
\end{equation}
Since it is very easy to observe that 

$$\left|
\frac{1}{\#}\int_1^{1+\#}
(e^{2\pi i \alpha\xi}-1) d\alpha
\right|$$

$$=\left(\left|\frac{1}{\#}\int_1^{1+\#}
(\cos(2\pi\alpha\xi)-1) d\alpha
\right|^2 + \left|
\frac{1}{\#}\int_1^{1+\#}
(\sin(2\pi\alpha\xi)) d\alpha
\right|^2\right)^{1/2}$$
is a number between $c_1$ and $2$ where $c_1>0$ is a universal constant, 
the claim follows.

Then, one can define the function $\Phi_I^{\text{left}}$ by

\begin{equation}
\widehat{\Phi_I^{\text{left}}}(\xi):=
\frac{\widehat{\Phi_I}(\xi)}{\widehat{\lambda_{\text{left}}}(\xi)}
\end{equation}
and this equality shows that $\Phi_I^{\text{left}}$ is indeed a bump adapted
to $I$. To verify the support condition, since
$\Phi_I^{\text{left}}\ast \lambda_{\text{left}} = \Phi_I$, it follows that

\begin{equation}\label{dif}
\frac{1}{\#}\int_1^{1+\#}
\left(\int_x^{x+\alpha}\Phi_I^{\text{left}}(t) dt\right) d\alpha = 0
\end{equation}
for every $x\geq b_I$. Differentiating (\ref{dif}) with respect to $x$, we 
obtain that

\begin{equation}\label{dif1}
\frac{1}{\#}\int_1^{1+\#}
\Phi_I^{\text{left}}(x+\alpha) d\alpha = \Phi_I^{\text{left}}(x)
\end{equation}
for every $x\geq b_I$. By iterating (\ref{dif1}) several times, we obtain
that for each $x\geq b_I$ and $n\in\N$ we have

\begin{equation}\label{iter}
\Phi_I^{\text{left}}(x)=...=
\frac{1}{\#}\int_1^{1+\#}...
\frac{1}{\#}\int_1^{1+\#}
\Phi_I^{\text{left}}(x+\alpha_1+...+\alpha_n) d\alpha_1... d\alpha_n.
\end{equation}
It is not difficult to see that the right hand side of (\ref{iter}) goes to
zero as $n$ goes to infinity since $\alpha_1,...,\alpha_n\in (1,1+\#)$
and $\Phi_I^{\text{left}}$ is a bump adapted to $I$ and this proves the
support condition. A similar argument works to treat the
$\lambda_{\text{right}}$, $\Phi_I^{\text{right}}$ case.
\end{proof} 
By ``rescaling'' Lemma \ref{lr}, one immediately obtains

\begin{corollary}\label{lrk}
Let $k\in\R$ and $I=[a_I, b_I]$ as in Lemma \ref{lr}. Denote by $I_k$ the 
interval $I_k:= [2^ka_I, 2^kb_I]$ and let $\Phi_{I_k}$ be a smooth function
supported on $I_k$. Define $\lambda^k_{\text{left}}$ and 
$\lambda^k_{\text{right}}$ by 
$\lambda^k_{\text{left}}(x):= 2^{-k}\lambda_{\text{left}}(2^{-k}x)$ and
$\lambda^k_{\text{right}}(x):= 2^{-k}\lambda_{\text{right}}(2^{-k}x)$.
Then, there exist bump functions $\Phi_{I_k}^{\text{left}}$ and
$\Phi_{I_k}^{\text{right}}$ adapted to $I_k$ so that

\begin{equation}\label{lk}
\Phi_{I_k}^{\text{left}}\ast \lambda^k_{\text{left}} = \Phi_{I_k}
\end{equation}
and

\begin{equation}\label{rk}
\Phi_{I_k}^{\text{right}}\ast \lambda^k_{\text{right}} = \Phi_{I_k}.
\end{equation}
Moreover, they have the additional properties that
$\supp  \Phi_{I_k}^{\text{left}}\subseteq (-\infty, 2^k b_I]$  and
$\supp  \Phi_{I_k}^{\text{right}}\subseteq [2^ka_I, +\infty )$.
\end{corollary} 
By ``modulating'' Corollary \ref{lrk}, one also obtains

\begin{corollary}\label{lrkn}
Let $k$, $I_k$, $\Phi_{I_k}$ be as in Corollary \ref{lrk} and let $n$ be a 
real number with $|n|\geq 1$.
Then, there exist bump functions $\Phi_{I_k}^{\text{left}, n}$ and
$\Phi_{I_k}^{\text{right}, n}$ adapted to $I_k$ so that

\begin{equation}\label{lkn}
n \left(\Phi_{I_k}^{\text{left},n}(\cdot) e^{2\pi i n 2^{-k}\cdot}\right)
\ast \lambda^k_{\text{left}}(x) = 
\Phi_{I_k}(x) e^{2\pi i n 2^{-k}x}
\end{equation}
and

\begin{equation}\label{rkn}
n \left(\Phi_{I_k}^{\text{right},n}(\cdot) e^{2\pi i n 2^{-k}\cdot}\right)
\ast \lambda^k_{\text{right}}(x) = 
\Phi_{I_k}(x) e^{2\pi i n 2^{-k}x}
\end{equation}
Moreover, they have the additional properties that
$\supp  \Phi_{I_k}^{\text{left}, n}\subseteq (-\infty, 2^k b_I]$  and also
$\supp  \Phi_{I_k}^{\text{right}, n}\subseteq [2^ka_I, +\infty )$.
\end{corollary} 

Coming back to our previous formula (\ref{desc}), let us now consider as there
a generic term of the form

\begin{equation}\label{generic}
D_{2^k}^{\infty}\widehat{\psi_j}(\xi) e^{2\pi i n \frac{9}{10} 2^{-k} \xi}
\end{equation}
for $n\in\Z$ and $j\in [-M, M]$ integer, with $j\neq 0$. To fix the situation,
assume also that $j>0$. By applying the ``right variant'' of 
Corollary \ref{lrkn}, we find a bump function $m^j_{k,n}$ adapted to the 
interval $[2^k(j-1), 2^k j]$, so that our term in (\ref{generic}) can be 
written as

$$n\left(
m^j_{k,n}(\cdot) e^{2\pi i n \frac{9}{10} 2^{-k}\cdot}
\ast \lambda^k_{\text{right}}\right)(\xi):=
n\left(
\widetilde{m^j_{k,n}}
\ast \lambda^k_{\text{right}}\right)(\xi)
$$

\begin{equation}\label{1}
=n 2^{-k} \frac{1}{\#}\int_1^{1+\#}
\left(\int_{\R}\chi_{[-\alpha 2^k, 0]}(t)\widetilde{m^j_{k,n}}(\xi+t) dt\right)
d\alpha.
\end{equation}
For a fixed $\alpha\in (1, 1+\#)$, since $m^j_{k,n}$ is a bump adapted to
$[2^k(j-1), 2^k j]$ and satisfying the support condition in 
Corollary \ref{lrkn}, it can be decomposed as

\begin{equation}\label{2}
m^j_{k,n} = \sum_{l=0}^{\infty} \frac{1}{(1+l)^{1000}} m^j_{k,n,l,\alpha}
\end{equation}
where $m^j_{k,n,l,\alpha}$ are smooth and uniformly bounded functions
supported inside $\frac{10}{9}[\alpha 2^k(j-1+l), \alpha 2^k(j+l)]$.
In particular, this allows us to decompose the term in (\ref{1}) as

\begin{equation}\label{3}
n \sum_{l=0}^{\infty} \frac{1}{(1+l)^{1000}}
2^{-k} \frac{1}{\#}\int_1^{1+\#}
\left(\int_{\R}\chi_{[-\alpha 2^k, 0]}(t)\widetilde{m^j_{k,n,l,\alpha}}(\xi+t) 
dt\right)
d\alpha. 
\end{equation}
Fix now $l$ and $\alpha$ and rewrite the corresponding inner integral in
(\ref{3}) in the form

\begin{equation}\label{4}
2^{-k}\int_{-\alpha 2^k}^0
\chi_{\omega_{\alpha, k}}(t)
\widetilde{m^n_{\omega^{j+l}_{\alpha, k}}}(\xi + t) dt
\end{equation}
where $\omega_{\alpha, k}$ and $\omega_{\alpha, k}^{j+l}$ denote the intervals
$[-\alpha 2^k, 0]$ and $[\alpha 2^k (j-1+l), \alpha 2^k (j+l)]$ 
respectively.

If $\omega$ is any other interval of the same length $\alpha 2^k$, then clearly
there exists a real number $t_0$ such that 
$\omega_{\alpha, k}^{j+l} + t_0 =\omega$ and we will denote by 
$\widetilde{m^n_{\omega}}$ the corresponding translated function
defined by the formula

\begin{equation}\label{5}
\widetilde{m^n_{\omega}}(t):= \widetilde{m^n_{\omega^{j+l}_{\alpha, k}}}(t-t_0)
\end{equation}
for any $t\in\R$. We shall also denote by $\D_k^{\alpha}$ the 
collection 
of all intervals of the form $[\alpha 2^k l', \alpha 2^k (l'+1)]$ for some
$l'\in\Z$. Then, we observe that (\ref{4}) is also equal to

\begin{equation}\label{6}
2^{-k}\int_{-\alpha 2^k}^0\left(\sum_{\omega\in \D^{\alpha}_{k}}
\chi_{\omega}(t)
\widetilde{m^n_{\omega^{l+j}}}(\xi + t)\right) dt
\end{equation}
where for a given $\omega\in \D^{\alpha}_{k}$, $\omega^{l+j}$ denotes
the interval $\omega^{l+j}:= \omega + (l+j)|\omega|$. On the other hand, 
by using
(\ref{5}), is is easy to observe that for a fixed $\xi\in\R$, the function

\begin{equation}\label{7}
t\rightarrow \sum_{\omega\in \D^{\alpha}_{k}}
\chi_{\omega}(t)
\widetilde{m^n_{\omega^{l+j}}}(\xi + t)
\end{equation}
is a periodic function of period $\alpha 2^k$. In particular, this implies that
the term in (\ref{6}) is also equal to

\begin{equation}\label{8}
2^{-k}\frac{1}{2L}
\int_{-L\alpha 2^k}^{L\alpha 2^k}
\left(\sum_{\omega\in \D^{\alpha}_{k}}
\chi_{\omega}(t)
\widetilde{m^n_{\omega^{l+j}}}(\xi + t)\right) dt
\end{equation}
for every $L\in\N$. We will denote from now on the limit as 
$L\rightarrow\infty$ of the expression in (\ref{8}) by

\begin{equation}\label{9}
\alpha \int{\!\!\!\!\!\! -}
\left(\sum_{\omega\in \D^{\alpha}_{k}}
\chi_{\omega}(t)
\widetilde{m^n_{\omega^{l+j}}}(\xi + t)\right) dt,
\end{equation}
and this essentially completes our desired decomposition.

If $j$ was negative in (\ref{generic}) then, we would have had to apply
the ``left variant'' of Corollary \ref{lrkn} instead.

To summarize, we managed to write our symbol $m(\xi_1, \xi_2)$ as

\begin{equation}\label{10}
m(\xi_1,\xi_2) = \sum_{\max(|j_1|, |j_2|)= M}
m^{j_1, j_2}(\xi_1, \xi_2)
\end{equation}
where each $m^{j_1, j_2}(\xi_1, \xi_2)$ is of the form

$$m^{j_1, j_2}(\xi_1, \xi_2) =$$

\begin{equation}\label{11}
\sum_{n_1, n_2\in\Z}
\sum_{l_1, l_2\in\N}
\frac{1}{(1+|n_1|)^{1000}}
\frac{1}{(1+|n_2|)^{1000}}
\frac{1}{(1+l_1)^{1000}}
\frac{1}{(1+l_2)^{1000}}
m^{j_1, j_2}_{n_1, n_2, l_1, l_2}(\xi_1, \xi_2)
\end{equation}
and where $m^{j_1, j_2}_{n_1, n_2, l_1, l_2}(\xi_1, \xi_2)$ are given by

$$m^{j_1, j_2}_{n_1, n_2, l_1, l_2}(\xi_1, \xi_2) =
\int_{\R} C^{j_1, j_2}_{k, n_1, n_2}\cdot$$

\begin{equation}\label{12}
\left(
\frac{1}{\#}
\int_1^{1+\#}\alpha_1
\left(
\int{\!\!\!\!\!\! -}
\left(\sum_{\omega_1\in \D^{\alpha_1}_{k}}
\chi_{\omega_1}(t_1)
\widetilde{m^{n_1}_{\omega_1^{l_1+j_1}}}(\xi_1 + t_1)\right) dt_1
\right) d\alpha_1
\right)\cdot
\end{equation}

$$\left(
\frac{1}{\#}
\int_1^{1+\#}\alpha_2
\left(
\int{\!\!\!\!\!\! -}
\left(\sum_{\omega_2\in \D^{\alpha_2}_{k}}
\chi_{\omega_2}(t_2)
\widetilde{m^{n_2}_{\omega_2^{l_2+j_2}}}(\xi_2 + t_2)\right) dt_2
\right) d\alpha_2
\right) dk
$$
while the constants $C^{j_1, j_2}_{k, n_1, n_2}$ satisfy

\begin{equation}\label{13}
|C^{j_1, j_2}_{k, n_1, n_2}|\lesssim 1
\end{equation}
uniformly in $k, n_1, n_2, j_1, j_2$.

The advantage of such a decomposition is that it is very well adapted to
arbitrary translations in the plane. More precisely, for $N=(N_1, N_2)\in\R^2$
one can write

$$\tau_N m^{j_1, j_2}_{n_1, n_2, l_1, l_2}(\xi_1, \xi_2) =
m^{j_1, j_2}_{n_1, n_2, l_1, l_2}(\xi_1 - N_1, \xi_2 - N_2)=
$$

$$\int_{\R} C^{j_1, j_2}_{k, n_1, n_2}\cdot$$

$$\left(
\frac{1}{\#}
\int_1^{1+\#}\alpha_1
\left(
\int{\!\!\!\!\!\! -}
\left(\sum_{\omega_1\in \D^{\alpha_1}_{k}}
\chi_{\omega_1}(t_1)
\widetilde{m^{n_1}_{\omega_1^{l_1+j_1}}}(\xi_1 - N_1 + t_1)\right) dt_1
\right) d\alpha_1
\right)\cdot$$

$$\left(
\frac{1}{\#}
\int_1^{1+\#}\alpha_2
\left(
\int{\!\!\!\!\!\! -}
\left(\sum_{\omega_2\in \D^{\alpha_2}_{k}}
\chi_{\omega_2}(t_2)
\widetilde{m^{n_2}_{\omega_2^{l_2+j_2}}}(\xi_2 - N_2 + t_2)\right) dt_2
\right) d\alpha_2
\right) dk = 
$$

$$\int_{\R} C^{j_1, j_2}_{k, n_1, n_2}\cdot$$

$$\left(
\frac{1}{\#}
\int_1^{1+\#}\alpha_1
\left(
\int{\!\!\!\!\!\! -}
\left(\sum_{\omega_1\in \D^{\alpha_1}_{k}}
\chi_{\omega_1}(t_1 + N_1)
\widetilde{m^{n_1}_{\omega_1^{l_1+j_1}}}(\xi_1 + t_1)\right) dt_1
\right) d\alpha_1
\right)\cdot$$

$$\left(
\frac{1}{\#}
\int_1^{1+\#}\alpha_2
\left(
\int{\!\!\!\!\!\! -}
\left(\sum_{\omega_2\in \D^{\alpha_2}_{k}}
\chi_{\omega_2}(t_2 + N_2)
\widetilde{m^{n_2}_{\omega_2^{l_2+j_2}}}(\xi_2 + t_2)\right) dt_2
\right) d\alpha_2
\right) dk = 
$$

$$\int_{\R} C^{j_1, j_2}_{k, n_1, n_2}\cdot$$

$$\left(
\frac{1}{\#}
\int_1^{1+\#}\alpha_1
\left(
\int{\!\!\!\!\!\! -}
\left(\sum_{\omega_1\in \D^{\alpha_1, t_1}_{k}}
\chi_{\omega_1}(N_1)
\widetilde{m^{n_1}_{\omega_1^{l_1+j_1}}}(\xi_1)\right) dt_1
\right) d\alpha_1
\right)\cdot$$

$$\left(
\frac{1}{\#}
\int_1^{1+\#}\alpha_2
\left(
\int{\!\!\!\!\!\! -}
\left(\sum_{\omega_2\in \D^{\alpha_2, t_2}_{k}}
\chi_{\omega_2}(N_2)
\widetilde{m^{n_2}_{\omega_2^{l_2+j_2}}}(\xi_2)\right) dt_2
\right) d\alpha_2
\right) dk = 
$$

\begin{equation}\label{14}
\frac{1}{\#}
\int_1^{1+\#}\alpha_1
\frac{1}{\#}
\int_1^{1+\#}\alpha_2
\int{\!\!\!\!\!\! -}
\int{\!\!\!\!\!\! -}
\int_{\R} C^{j_1, j_2}_{k, n_1, n_2}\cdot
\end{equation}

$$
\left(
\sum_{\omega_1\in \D^{\alpha_1, t_1}_{k}}
\sum_{\omega_2\in \D^{\alpha_2, t_2}_{k}}
\chi_{\omega_1}(N_1)
\chi_{\omega_2}(N_2)
\widetilde{m^{n_1}_{\omega_1^{l_1+j_1}}}(\xi_1)
\widetilde{m^{n_2}_{\omega_2^{l_2+j_2}}}(\xi_2)
\right) dk
 d t_1 d t_2 d\alpha_1 d\alpha_2 =
$$

\begin{equation}\label{14'}
\frac{1}{\#}
\int_1^{1+\#}\alpha_1
\frac{1}{\#}
\int_1^{1+\#}\alpha_2
\int{\!\!\!\!\!\! -}
\int{\!\!\!\!\!\! -}
\int_0^1 \sum_{k\in\Z} C^{j_1, j_2}_{k+\kappa, n_1, n_2}\cdot
\end{equation}

$$
\left(
\sum_{\omega_1\in \D^{\alpha_1, t_1}_{k+\kappa}}
\sum_{\omega_2\in \D^{\alpha_2, t_2}_{k+\kappa}}
\chi_{\omega_1}(N_1)
\chi_{\omega_2}(N_2)
\widetilde{m^{n_1}_{\omega_1^{l_1+j_1}}}(\xi_1)
\widetilde{m^{n_2}_{\omega_2^{l_2+j_2}}}(\xi_2)
\right) d\kappa
 d t_1 d t_2 d\alpha_1 d\alpha_2,
$$
where in general, $\D^{\alpha, t}_{k}$ denotes the set of all intervals of 
the form
$\omega - t$ with $\omega\in \D^{\alpha}_{k}$.

As a consequence of our decomposition (\ref{10}) - (\ref{14'}) 
our maximal operator $C_m (f_1, f_2)(x)$ can be estimated by

\begin{equation}\label{15}
C_m (f_1, f_2)(x)\lesssim
\sum_{\max(|j_1|, |j_2|) = M}
C^{j_1, j_2}(f_1, f_2)(x)
\end{equation}
while each $C^{j_1, j_2}(f_1, f_2)(x)$ can be majorized by

\begin{equation}\label{16}
C^{j_1, j_2}(f_1, f_2)(x)\lesssim
\end{equation}

$$
\sum_{n_1, n_2\in\Z}
\sum_{l_1, l_2\in\N}
\frac{1}{(1+|n_1|)^{1000}}
\frac{1}{(1+|n_2|)^{1000}}
\frac{1}{(1+l_1)^{1000}}
\frac{1}{(1+l_2)^{1000}}
C^{j_1, j_2}_{n_1, n_2, l_1, l_2}(f_1, f_2)(x)
$$
where $C^{j_1, j_2}_{n_1, n_2, l_1, l_2}(f_1, f_2)(x)$ is the maximal operator
defined by the formula (\ref{carlesonm-def}) for the case of the symbol
$m^{j_1, j_2}_{n_1, n_2, l_1, l_2}(\xi_1, \xi_2)$.

Fix now $(p_1, p_2, p)$ as in Theorem \ref{main}. Clearly, to prove our main
theorem, it is enough to prove the inequality

\begin{equation}\label{17}
\|C^{j_1, j_2}_{n_1, n_2, l_1, l_2}(f_1, f_2)\|_p
\lesssim
(1+|n_1|)^{10}
(1+|n_2|)^{10}
(1+l_1)^{10}
(1+l_2)^{10}
\|f_1\|_{p_1}
\|f_2\|_{p_2}
\end{equation}
for every $f_1\in L^{p_1}$ and $f_2\in L^{p_2}$. 
To avoid unnecessary technical complications,
we are only going to prove 
that the operator $C^{j_1, j_2}_{0,0,0,0}$ satisfies the desired estimates,
but it will be clear from the proof we shall present that the same arguments 
give the general inequality (\ref{17}). 

We will therefore concentrate our attention on the operator
$C^{j_1, j_2}_{0,0,0,0}$, from now on.

\section{Restricted type estimates}

The purpose of this section is to review the interpolation theory from
\cite{cct} which allows us to reduce the estimates in Theorem \ref{main}
to certain ``restricted type estimates''. Roughly speaking, we will see that 
it is enough to prove our desired estimates in the particular case when all the
functions involved are characteristic functions of measurable sets.

Consider now $j_1, j_2$ integers so that $\max(|j_1|, |j_2|) = M$ and to fix
the case assume from now on that both $j_1$ and $j_2$ are positive
(all the other cases can be treated in the same way).

Then, we linearize our operator $C^{j_1, j_2}_{0,0,0,0}$ in (\ref{17}) as

\begin{equation}\label{18}
C^{j_1, j_2}_{0,0,0,0}(f_1, f_2)(x) =
\int_{\R^2} \tau_{N(x)} m^{j_1, j_2}_{0,0,0,0}(\xi_1, \xi_2)
\widehat{f_1}(\xi_1)
\widehat{f_2}(\xi_2)
e^{2\pi i x (\xi_1 + \xi_2)} d\xi_1 d\xi_2.
\end{equation}
To prove the $L^p$ estimates on $C^{j_1, j_2}_{0,0,0,0}$, it is convenient
to use duality and introduce the trilinear form $\Lambda^{j_1, j_2}$
via the formula

\begin{equation}\label{form}
\Lambda^{j_1, j_2}(f_1, f_2, f_3):= \int_{\R}
C^{j_1, j_2}_{0,0,0,0}(f_1, f_2)(x) f_3(x) dx.
\end{equation}
Then, the statement that $C^{j_1, j_2}_{0,0,0,0}$ is bounded from
$L^{p_1}\times L^{p_2}\rightarrow L^p$ is equivalent to $\Lambda^{j_1, j_2}$
being bounded on $L^{p_1}\times L^{p_2}\times L^{p'}$ if $1\leq p <\infty$.
For $p<1$ this duality does no longer hold, however the interpolation arguments
in \cite{cct} will allow us to replace it with certain restricted type 
estimates on $\Lambda^{j_1, j_2}$. As in \cite{cct}, we find it more
convenient to work with the quantities $\alpha_1=\frac{1}{p_1}$,
$\alpha_2=\frac{1}{p_2}$ and $\alpha_3 = \frac{1}{p'}$ where $p_1, p_2, p$
stand for the exponents of the spaces $L^{p_1}$, $L^{p_2}$ and $L^p$.

We recall now the following definitions which have been introduced in 
\cite{cct}.

\begin{definition}
A tuple $\alpha=(\alpha_1,\alpha_2,\alpha_3)$ is called admissible, if

\[-\infty <\alpha_i <1\]
for all $1\leq i\leq 3$,

\[\sum_{i=1}^3\alpha_i=1\]
and there is at most one index $j$ such that $\alpha_j<0$. We call
an index $i$ good if $\alpha_i\geq 0$, and we call it bad if
$\alpha_i<0$. A good tuple is an admissible tuple without bad index, a
bad tuple is an admissible tuple with only one bad index.
\end{definition}

\begin{definition}
Let $E$, $E'$ be sets of finite measure. We say that $E'$ is a major 
subset of $E$ if $E'\subseteq E$ and $|E'|\geq\frac{1}{2}|E|$.
\end{definition} 

\begin{definition}
If $E$ is a set of finite measure, we denote by $X(E)$ the space of
all measurable complex-valued functions $f$ supported on $E$ and such that 
$\|f\|_{\infty}\leq
1$.
\end{definition}

\begin{definition}
If $\alpha=(\alpha_1,\alpha_2,\alpha_3)$ is an admissible bad tuple
with bad index $j$, we say
that our $3$-linear form $\Lambda^{j_1, j_2}$ is of restricted type $\alpha$
if for every sequence $E_1, E_2, E_3$ of subsets of $\R$ with
finite measure, there exists a major subset $E'_j$ of $E_j$ such that

\[|\Lambda^{j_1, j_2}(f_1, f_2, f_3)|\lesssim |E|^{\alpha}\]
for all functions $f_i\in X(E'_i)$, $i=1,2,3$, where we adopt the
convention $E'_i =E_i$ for good indices $i$, and $|E|^{\alpha}$
is a shorthand for

\[|E|^{\alpha}=|E_1|^{\alpha_1}|E_2|^{\alpha_2}
|E_3|^{\alpha_3}.\]
If $\alpha=(\alpha_1,\alpha_2,\alpha_3)$ is an admissible good tuple, 
we say
that the $3$-linear form $\Lambda^{j_1,j_2}$ is of restricted type $\alpha$ if
there exists $j$ such that
for every sequence $E_1, E_2, E_3$ of subsets of $\R$ with
finite measure, there exists a major subset $E'_j$ of $E_j$ such that

\[|\Lambda^{j_1, j_2}(f_1, f_2, f_3)|\lesssim |E|^{\alpha}\]
for all functions $f_i\in X(E'_i)$, $i=1,2,3$, where this time we adopt the
convention $E'_i =E_i$ for the indices $i\neq j$.

\end{definition}

Let us consider now the $2$-dimensional affine hyperplane
\[
S:=\{(\alpha_1,\alpha_2,\alpha_3)\in\R^3\,
|\,\alpha_1 + \alpha_2 + \alpha_3=1\}.
\]

The points $A_1, A_2, A_3$ and $B_1, B_2, B_3$ belong to $S$ and have the 
following coordinates:

\[
\begin{array}{llll}
A_1:(1, 1, -1)  &  A_2:(-1, 1, 1)  &  
A_3:(1, -1, 1)   \\  
\ &\ &\ & \ \\
B_1:(0, 1, 0)  &  B_2:(0, 0, 1) & 
B_3:(1, 0, 0).   \\  
\end{array}
\]

\setlength{\unitlength}{1.1mm}
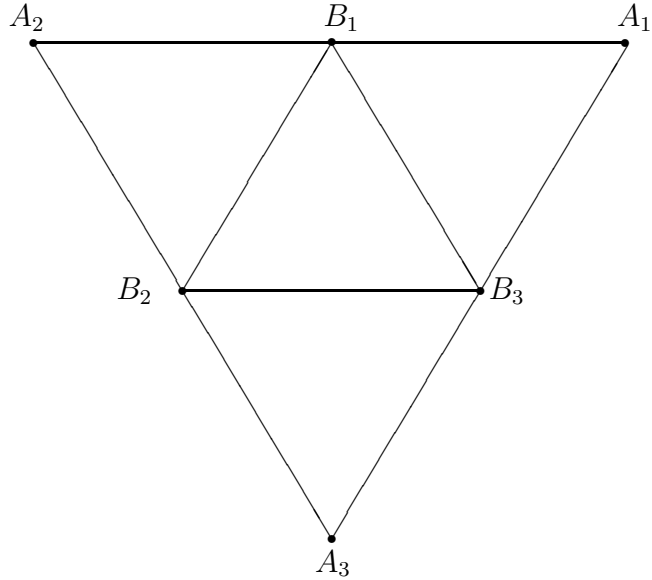
\begin{figure}\label{Triangle}
\caption{Triangle}
\begin{picture}(72,78)
\put(46,60){\circle*{0.8}}
\put(28,30){\circle*{0.8}}
\put(64,30){\circle*{0.8}}
\put(81.5,60){\circle*{0.8}}
\put(45,62){$B_1$}
\put(80.5,62){$A_1$}
\put(20,29){$B_2$}
\put(65,29){$B_3$}
\put(7,62){$A_2$}
\put(10,60){\circle*{0.8}}
\put(44,-4){$A_3$}
\put(46,0){\circle*{0.8}}



\put(10,60){\line(1,0){71.4}}
\put(28.4,30){\line(1,0){35.2}}
\put(10,60){\line(3,-5){35.8}}
\put(46.2,59.67){\line(3,-5){17.6}}
\put(28.2,30.33){\line(3,5){17.6}}
\put(46,0){\line(3,5){35.8}}




\end{picture}
\end{figure}

The following restricted type estimates will be proved directly.

\begin{theorem}\label{main'}
For every vertex $A_1, A_2, A_3$ there exist admissible tuples $\alpha$
arbitrarily close to $A_1, A_2, A_3$ such that the trilinear form
$\Lambda^{j_1, j_2}$ is of restricted type $\alpha$.
\end{theorem}
By interpolation of these restricted type estimates we obtain as in \cite{cct}
the following corollary.

\begin{corollary}
Let $\alpha$ be an admissible tuple inside the interior of the triangle
$[A_1 A_2 A_3]$. Then, $\Lambda^{j_1, j_2}$ is of restricted type $\alpha$.
\end{corollary}
Since one observes that $p_1, p_2, p$ satisfy the hypothesis of Theorem
\ref{main} if and only if 
$(\frac{1}{p_1}, \frac{1}{p_2}, \frac{1}{p'})\in [A_1 B_1 B_2 B_3]$, it only
remains to convert these restricted type estimates into strong type estimates.

To do this, as in \cite{cct}, \cite{mtt:fourierbiest}, \cite{mtt:bicarleson},
one just has to apply exactly as in \cite{cct}, the multilinear 
Marcinkiewicz interpolation theorem in the case of good tuples and Lemma 3.1
in \cite{cct} in the case of bad tuples.

In other words, we have proved that our main Theorem \ref{main} can be 
reduced to Theorem \ref{main'}. It thus remains to only prove 
Theorem \ref{main'}.

\section{A discretized model}

The main task of the present section is to reduce Theorem \ref{main'} to
a discretized variant of it. We first introduce more notations.

Let $j_1, j_2$ as before, $\alpha_1, \alpha_2\in [1,1+\#]$, $t_1, t_2\in\R$
and $\kappa\in [0,1]$. For any $k\in\Z$ we denote by
$\vec{\overline{\P}}^{j_1, j_2}_{\alpha_1, \alpha_2, t_1, t_2, \kappa, k}$
the set of all $4$ - tuples of the form $\vec{P} = (P_1, P_2, P_3, P_4)$
where $P_i$ are tiles defined by $P_i:= I_{\vec{P}}\times \omega_{P_i}$
with $I_{\vec{P}}\in \D^{0, 0}_{-k}$, 
$\omega_{P_i}\in \D^{\alpha_1, t_1}_{k+\kappa}$ for $i=1,3$ and
$\omega_{P_i}\in \D^{\alpha_2, t_2}_{k+\kappa}$ for $i=2,4$. 
Moreover, they also have the important property that 
$\omega_{P_3} = \omega_{P_1} - j_1|\omega_{P_1}|$ and
$\omega_{P_4} = \omega_{P_2} - j_2|\omega_{P_2}|$. Then, we denote by
$\vec{\overline{\P}}^{j_1, j_2}_{\alpha_1, \alpha_2, t_1, t_2, \kappa}$
the set

\begin{equation}
\vec{\overline{\P}}^{j_1, j_2}_{\alpha_1, \alpha_2, t_1, t_2, \kappa}:=
\bigcup_{k\in\Z}
\vec{\overline{\P}}^{j_1, j_2}_{\alpha_1, \alpha_2, t_1, t_2, \kappa, k}
\end{equation}
which we name the set of all vector-tiles of dimension $4$ associated to the
parameters $\alpha_1, \alpha_2, t_1, t_2, \kappa$. Note that our 
generic collection of intervals 
$\cup_{k\in\Z}\D^{\alpha, t}_{k+\kappa}$ forms what is usually 
called a dyadic grid. In particular, if 
$I, J\in \cup_{k\in\Z}\D^{\alpha, t}_{k+\kappa} $ have the property that
$I\cap J\neq\emptyset$, then either $I\subseteq J$ or $J\subseteq I$. We will
freely use this facts throughout the paper.

The discrete version of Theorem \ref{main'} is the following.

\begin{theorem}\label{discrete}
Let $j_1, j_2, \alpha_1, \alpha_2, t_1, t_2, \kappa$ be as before and let
$\vec{\P}^{j_1, j_2}_{\alpha_1, \alpha_2, t_1, t_2, \kappa}\subseteq
\vec{\overline{\P}}^{j_1, j_2}_{\alpha_1, \alpha_2, t_1, t_2, \kappa}$ be a 
finite collection of vector-tiles of dimension $4$ associated to them. For each
$\vec{P}\in \vec{\P}^{j_1, j_2}_{\alpha_1, \alpha_2, t_1, t_2, \kappa}$ and
$i=1,2$ let $\Phi_{P_i} = \Phi_{P_i, i}$ be an $L^2$ normalized  bump
function adapted to $I_{\vec{P}}$ with 
$\supp \widehat{\Phi_{P_i} }\subseteq \frac{10}{9}\omega_{P_i}$ and let also
$\Phi_{\vec{P}}$ be an $L^2$ normalized bump function adapted to $I_{\vec{P}}$
with $\supp \widehat{\Phi_{\vec{P}}}\subseteq 
\frac{9}{8}(-\omega_{P_1}-\omega_{P_2})$. Define the trilinear form
$\Lambda_{\vec{\P}^{j_1, j_2}_{\alpha_1, \alpha_2, t_1, t_2, \kappa} }$ by

\begin{equation}\label{discreteform}
\Lambda_{\vec{\P}^{j_1, j_2}_{\alpha_1, \alpha_2, t_1, t_2, \kappa} }
(f_1, f_2, f_3) =
\end{equation}

$$\sum_{\vec{P}\in \vec{\P}^{j_1, j_2}_{\alpha_1, \alpha_2, t_1, t_2, \kappa}}
\frac{1}{|I_{\vec{P}}|^{1/2}}
\langle f_1, \Phi_{P_1}\rangle
\langle f_2, \Phi_{P_2}\rangle
\langle f_3, \Phi_{\vec{P}}
\chi_{\{x: N_1(x)\in\omega_{P_3}\}}
\chi_{\{x: N_2(x)\in\omega_{P_4}\}}\rangle.
$$
Then, for every vertex $A_j$, $j=1,2,3$ there exist admissible tuples $\alpha$
arbitrarily close to $A_j$ such that the form 
$\Lambda_{\vec{\P}^{j_1, j_2}_{\alpha_1, \alpha_2, t_1, t_2, \kappa} }$
is of restricted type $\alpha$ uniformly with respect to the parameters
$\alpha_1, \alpha_2, t_1, t_2, \kappa$, 
$\vec{\P}^{j_1, j_2}_{\alpha_1, \alpha_2, t_1, t_2, \kappa}$,
$\Phi_{P_i}$, $\Phi_{\vec{P}}$. Moreover, in the case when $\alpha$
has bad index $j$, the restricted type is uniform in the sense that the major
subset $E'_j$ can be chosen independently of all the parameters just mentioned.
\end{theorem}
In the rest of this section we show how Theorem \ref{main'} can be deduced
from Theorem \ref{discrete}. To realize this, one just has to calculate
the form $\Lambda^{j_1, j_2}$ in (\ref{form}) carefully.

First, we fix $N=(N_1, N_2)\in\R^2$ and look at the trilinear form associated
to the bilinear operator $T_{\tau_N m^{j_1, j_2}_{0,0,0,0}}$. It is given by
the formula

\begin{equation}
\int_{\R^3}\delta(\xi_1, \xi_2, \xi_3)
\tau_N m^{j_1, j_2}_{0,0,0,0}(\xi_1, \xi_2)
\widehat{f_1}(\xi_1)
\widehat{f_2}(\xi_2)
\widehat{f_3}(\xi_3)
d\xi_1 d\xi_2 d\xi_3
\end{equation}
where $\delta$ is the Dirac delta function. We see from (\ref{14'}) that the
symbol $\tau_N m^{j_1, j_2}_{0,0,0,0}(\xi_1, \xi_2)$ is an average over
parameters $\alpha_1, \alpha_2, t_1, t_2, \kappa$ of sums over $k\in\Z$,
$\omega_1\in \D^{\alpha_1, t_1}_{k+\kappa}$ and
$\omega_2\in \D^{\alpha_2, t_2}_{k+\kappa}$. We fix all these parameters
and fix also $k, \omega_1, \omega_2$. The corresponding multiplier, is given
by the expression

\begin{equation}
\chi_{\omega_1}(N_1)
\chi_{\omega_2}(N_2)
m_{\omega_1^{j_1}}(\xi_1)
m_{\omega_2^{j_2}}(\xi_2)
\end{equation}
and as a consequence, the trilinear form associated to it has the formula

\begin{equation}\label{18}
\int_{\R^3}
\delta(\xi_1, \xi_2, \xi_3)
\widehat{f_1}(\xi_1)
m_{\omega_1^{j_1}}(\xi_1)
\widehat{f_2}(\xi_2)
m_{\omega_2^{j_2}}(\xi_2)
\widehat{f_3}(\xi_3)
\chi_{\omega_1}(N_1)
\chi_{\omega_2}(N_2)
d\xi_1 d\xi_2 d\xi_3.
\end{equation}
Pick now $m_{-\omega_1^{j_1}-\omega_2^{j_2}}$ a smooth function supported
inside $\frac{9}{8}(-\omega_1^{j_1}-\omega_2^{j_2} )$ and which equals $1$
on $-\frac{10}{9}\omega_1^{j_1}-\frac{10}{9}\omega_2^{j_2}$. Then, formula
(\ref{18}) can also be written as

\begin{equation}
\chi_{\omega_1}(N_1)
\chi_{\omega_2}(N_2)
\int_{\R^3}
\delta(\xi_1, \xi_2, \xi_3)
\widehat{f_1\ast\check{m_{\omega_1^{j_1}}}}
\widehat{f_2\ast\check{m_{\omega_2^{j_2}}}}
\widehat{f_3\ast\check{m_{-\omega_1^{j_1}-\omega_2^{j_2}}}}
d\xi_1 d\xi_2 d\xi_3.
\end{equation}
By using Plancherel, this is equal to

\begin{equation}
\chi_{\omega_1}(N_1)
\chi_{\omega_2}(N_2)
\int_{\R}
\left(f_1\ast\check{m_{\omega_1^{j_1}}}\right)(x)
\left(f_2\ast\check{m_{\omega_2^{j_2}}}\right)(x)
\left(f_3\ast\check{m_{-\omega_1^{j_1}-\omega_2^{j_2}}}\right)(x) dx=
\end{equation}

$$2^{3/2 k}\int_{\R}
\langle f_1, \Phi^1_{\omega, x}\rangle
\langle f_2, \Phi^2_{\omega, x}\rangle
\langle f_3, \Phi^3_{\omega, x}\rangle dx
$$
where

$$\Phi^i_{\omega, x}(y):= 
2^{-k/2}\overline{\check{m_{\omega_i^{j_i}}}(x-y) }
$$
for $i=1,2$ and

$$\Phi^3_{\omega, x}(y):= 
2^{-k/2}\overline{\check{m_{-\omega_1^{j_1}-\omega_2^{j_2} }}(x-y) }.
$$
We can rewrite this as

\begin{equation}
\int_0^1
\sum_{\vec{P}}
\frac{1}{|I_{\vec{P}}|^{1/2}}
\langle f_1, \Phi_{P_1, t, 1}\rangle
\langle f_2, \Phi_{P_2, t, 2}\rangle
\langle f_3, \Phi_{P_3, t, 3}\rangle
\chi_{\omega_{P_3}}(N_1)
\chi_{\omega_{P_4}}(N_2)
dt
\end{equation}
where $\vec{P}$ ranges over all vector-tiles of dimension $4$ having the 
property that $\omega_{P_1}= \omega_1^{j_1}$,
$\omega_{P_2}= \omega_2^{j_2}$, $\omega_{P_3}=\omega_1$,
$\omega_{P_4}=\omega_2$, $I_{\vec{P}}$ is in $\D_{-k}^{0,0}$,
$\Phi_{P_j, t, j}$ is the function

$$\Phi_{P_j, t, j}:= \Phi^j_{\omega_{P_j}, x_{\vec{P}}+t|I_{\vec{P}}|}$$
and $x_{\vec{P}}$ is the center of $I_{\vec{P}}$.
As a consequence of these computations, it follows that the bilinear
operator $T_{\tau_N m^{j_1, j_2}_{0,0,0,0}}(f_1, f_2)$ can be written as 
an average over parameters $\alpha_1, \alpha_2, t_1, t_2, \kappa$ of
expressions of the form

\begin{equation}
\sum_{\vec{P}\in
\vec{\overline{\P}}^{j_1, j_2}_{\alpha_1, \alpha_2, t_1, t_2, \kappa}}
\frac{1}{|I_{\vec{P}}|^{1/2}}
\langle f_1, \Phi_{P_1, t, 1}\rangle
\langle f_2, \Phi_{P_2, t, 2}\rangle
\overline{\Phi_{P_3, t, 3} }
\chi_{\omega_{P_3}}(N_1)
\chi_{\omega_{P_4}}(N_2)
\end{equation}
and in particular this means that the linearized operator 
$C^{j_1, j_2}_{0,0,0,0}$ can be written as an average over the same parameters
of expressions of the form

\begin{equation}
\sum_{\vec{P}\in
\vec{\overline{\P}}^{j_1, j_2}_{\alpha_1, \alpha_2, t_1, t_2, \kappa}}
\frac{1}{|I_{\vec{P}}|^{1/2}}
\langle f_1, \Phi_{P_1, t, 1}\rangle
\langle f_2, \Phi_{P_2, t, 2}\rangle
\overline{\Phi_{P_3, t, 3} }
\chi_{\omega_{P_3}}(N_1(x))
\chi_{\omega_{P_4}}(N_2(x)).
\end{equation}

Now it is clear that to prove our claim one should simply integrate
the conclusion of Theorem \ref{discrete} over all the above parameters
using the uniformity assumptions of that theorem. The finiteness condition
on $\vec{\P}^{j_1, j_2}_{\alpha_1, \alpha_2, t_1, t_2, \kappa}$
can be easily removed by a standard limiting argument.

\section{Trees and vector-trees}

The standard approach to prove our desired estimates on the form
$\Lambda_{\vec{\P}^{j_1, j_2}_{\alpha_1, \alpha_2, t_1, t_1, \kappa}}$
is to organize our fixed collection of vector-tiles of dimension $4$
into vector-trees. To define them rigorously we need to recall some ordering 
relations between tiles as in \cite{fefferman}, \cite{laceyt1},
\cite{cct}.

Let $\alpha\in (1,1+\#)$, $t\in\R$ and $\kappa\in [0,1]$. We will denote
by $\overline{\P}_{\alpha, t}$ the set defined by

\begin{equation}
\overline{\P}_{\alpha, t}:=
\bigcup_{k\in\Z}
\D^{0,0}_{-k}\times \D^{\alpha, t}_{k+\kappa}
\end{equation}
of all tiles $P= I_P\times \omega_P$.

\begin{definition}\label{order}
Let $P$, $P'$ be tiles in $\overline{\P}_{\alpha, t}$. We define the ordering
$<$ and write $P' < P$ if $I_{P'}\subsetneq I_P$ and 
$3\omega_P\subseteq 3\omega_{P'}$ and $P'\leq P$ if $P' < P$ or $P'=P$. We
also write $P'\lesssim P$ if $I_{P'}\subsetneq I_P$ and
$3M\omega_P\subseteq 3M\omega_{P'}$. Finally, we also write $P'\lesssim' P$
if $P'\lesssim P$ and $P'\not \leq P$.

We will sometime also use the classical C.Fefferman's ordering \cite{fefferman}
and write $P'<^c P$ if $I_{P'}\subsetneq I_P$ and 
$\omega_P\subseteq \omega_{P'}$.
\end{definition}

\begin{definition}\label{tree}
Let $\P_{\alpha, t}\subseteq \overline{\P}_{\alpha, t}$ be an arbitrary 
collection of tiles and let $P_0\in \overline{\P}_{\alpha, t}$ be a fixed tile.
A collection $T\subseteq \P_{\alpha, t}$ is called a tree with top $P_0$
if and only if

\begin{equation}
P\leq P_0
\end{equation}
for all $P\in T$. We write $I_T$ for $I_{P_0}$. Note that the tree $T$ does
not necessarily contain its top $P_0$.

Similarly, a collection $T\subseteq \P_{\alpha, t}$ is called a lacunary tree
with top $P_0$ if and only if 

\begin{equation}
P\lesssim' P_0
\end{equation}
for all $P\in T$.
\end{definition}

\begin{definition}\label{disjoint}
Two trees $T$ and $T'$ are said to be strongly disjoint if and only if

(a) $P\neq P'$ for all $P\in T$ and $P'\in\T'$.

(b) Whenever $P\in T$ and $P'\in T'$ are such that 
$3\omega_P\cap 3\omega_{P'}\neq\emptyset$ then, one has 
$I_{P'}\cap I_T=\emptyset$ and similarly with $T$ and $T'$ reversed.
\end{definition}

\begin{definition}\label{sparse}
An arbitrary collection of tiles $\P_{\alpha, t}$ is called sparse
if and only if for any two tiles $P, P'\in \P_{\alpha, t}$ we have
$|\omega_P|< |\omega_{P'}|$ implies $2M|\omega_P|< |\omega_{P'}|$ and
$|\omega_P|=|\omega_{P'}|$ implies $M\omega_P\cap M\omega_{P'}=\emptyset$.
\end{definition}

It is very easy to observe that any given collection of tiles can be written 
as an $O(M^2)$ disjoint union of sparse collections of tiles. From now on,
we will assume throughout the paper that all our collections of tiles
are sparse.

\begin{definition}\label{vectortree4}
Let $j_1, j_2, \alpha_1, \alpha_2, t_1, t_2, \kappa$ be as before and let
$\vec{\P}^{j_1, j_2}_{\alpha_1, \alpha_2, t_1, t_2, \kappa}\subseteq
\vec{\overline{\P}}^{j_1, j_2}_{\alpha_1, \alpha_2, t_1, t_2, \kappa}$ be an 
arbitrary collection of vector-tiles of dimension $4$ associated to them.
A collection $\vec{T}\subseteq
\vec{\P}^{j_1, j_2}_{\alpha_1, \alpha_2, t_1, t_2, \kappa}$ is called a 
vector-tree of dimension $4$ if and only if for each $1\leq i\leq 4$ the
projected collection

\begin{equation}
T_i:= \{ P_i : \vec{P}=(P_1, P_2, P_3, P_4)\in \vec{T} \}
\end{equation}
is a tree of tiles (lacunary or not).
\end{definition}
From now on, for the simplicity of our notation we will omit the indices
$j_1, j_2, \alpha_1, \alpha_2, t_1, t_2, \kappa$ and simply write
$\vec{\P}\subseteq\vec{\overline{\P}}$ instead of
$\vec{\P}^{j_1, j_2}_{\alpha_1, \alpha_2, t_1, t_2, \kappa}\subseteq
\vec{\overline{\P}}^{j_1, j_2}_{\alpha_1, \alpha_2, t_1, t_2, \kappa}$.

\begin{definition}\label{vectortree2}
Let $\vec{\P}\subseteq\vec{\overline{\P}}$ as before and denote by
$\vec{\P}_{3,4}\subseteq\vec{\overline{\P}}_{3,4}$ the sets of projected
vector-tiles of dimension $2$ defined by

$$\vec{\P}_{3,4}:= \{ (P_3, P_4) : \vec{P}=(P_1, P_2, P_3, P_4)\in\vec{\P}\}$$
and

$$\vec{\overline{\P}}_{3,4}:= \{ (P_3, P_4) : \vec{P}=(P_1, P_2, P_3, P_4)
\in\vec{\overline{\P}}\}.$$
A collection $\vec{T}\subseteq \vec{\P}_{3,4}$ is called a vector-tree 
of dimension $2$ if and only if for each $i=3,4$ the projected collection

$$T_i:= \{ P_i : (P_3, P_4)\in\T \}$$
is a tree of tiles (lacunary or not).
\end{definition}

\section{Sizes and energies}

The standard way to estimate our trilinear form $\Lambda_{\vec{\P}}$, 
is to do so by introducing some ``sizes'' and ``energies'' well adapted to our
given collection of vector-tiles. The first, have been considered in
\cite{mtt:fourierbiest}.

\begin{definition}\label{se12}
Let $\vec{\P}\subseteq\vec{\overline{\P}}$ be as usual and $1\leq i\leq 2$. 
We denote
by $\P_i$ the set of all distinct tiles $P_i$ where 
$\vec{P}=(P_1, P_2, P_3, P_4)\subseteq\vec{\P}$. 
Then, for every function $f_i$ we define the size of it by

\begin{equation}
\size_{i,\P_i}(f_i):= \sup_{T\subseteq \P_i}
\left(\frac{1}{|I_T|} \sum_{P\in T}
|\langle f_i, \Phi_{P_i}\rangle|^2
\right)^{1/2}
\end{equation}
where the suppremum is taken over all lacunary trees $T\subseteq \P_i$.

We also define the energy of the function $f_i$ by

\begin{equation}
\energy_{i,\P_i}(f_i):=\sup_{n\in\Z} \sup_{\bf{D}} 2^n
\left(\sum_{T\in\bf{D}} |I_T|\right)^{1/2}
\end{equation}
where $\bf{D}$ ranges over all collections of strongly disjoint 
lacunary trees in
$\P_i$ such that

$$
\left(\sum_{P\in T}
|\langle f_i, \Phi_{P_i}\rangle|^2
\right)^{1/2}\geq 2^n |I_T|^{1/2}
$$
for all $T\in\bf{D}$ and 

$$
\left(\sum_{P\in T'}
|\langle f_i, \Phi_{P_i}\rangle|^2
\right)^{1/2}\leq 2^{n+1} |I_T|^{1/2}
$$
for all lacunary subtrees $T'\subseteq T$.
\end{definition} 

Clearly, the sizes are phase-space variants of the $BMO$ norm of $f_i$, while
the energies are phase-space variants of the $L^2$ norm of $f_i$, for $i=1,2$.

The following John-Nirenberg type inequality is also true \cite{cct}.

\begin{lemma}\label{jn}
If $\vec{\P}$ and $\P_i$, $i=1,2$ are as before, then

$$\size_{i,\P_i}(f_i) \sim
\sup_{T\subseteq \P_i}
\frac{1}{|I_T|}
\left\|\left(
\sum_{P\in T}
\frac{|\langle f_i, \Phi_i\rangle|^2}{|I_{P_i}|}\chi_{I_{P_i}}\right)^{1/2}
\right\|_{1,\infty}
$$
where again, $T$ ranges over all lacunary trees in $\P_i$.
\end{lemma}

We now need sizes and energies to take care of our third function $f_3$.
They are defined as follows.

\begin{definition}\label{se34}
Let $\vec{\P}\subseteq\vec{\overline{\P}}$ be an arbitrary set of 
vector-tiles of dimension $4$ and 
$\vec{\P}_{3,4}\subseteq \vec{\overline{\P}}_{3,4}$ be the corresponding
set of vector-tiles of dimension $2$. We define the size of the function
$f_3$ by

\begin{equation}
\size_{3,\vec{\P}_{3,4}}(f_3):=
\sup_{\vec{P}\in \vec{\P}_{3,4}}
\sup_{\vec{P'}\in \vec{\overline{\P}}_{3,4}(\vec{P})}
\frac{1}{|I_{\vec{P'}}|}
\int_{\R}
|f_3(x)| \widetilde{\chi}_{I_{\vec{P'}}}^C(x)\cdot
\end{equation}

$$\cdot
\chi_{\{ x : N_1(x)\in \cup_{j=0}^{j_1}(\omega_{P'_3}+j|\omega_{P'_3}|) \}}
\cdot
\chi_{\{ x : N_2(x)\in \cup_{j=0}^{j_2}(\omega_{P'_4}+j|\omega_{P'_4}|) \}}
dx
$$
where $\vec{\overline{\P}}_{3,4}(\vec{P})$ is the set of all 
$\vec{P'}\in \vec{\overline{\P}}_{3,4}$ having the property that
$I_{\vec{P}}\subsetneq I_{\vec{P'}}$ and also that

$$\bigcup_{j=0}^{j_1}(\omega_{P'_3}+j|\omega_{P'_3}|)\subseteq
\bigcup_{j=0}^{j_1}(\omega_{P_3}+j|\omega_{P_3}|)$$
and

$$\bigcup_{j=0}^{j_2}(\omega_{P'_4}+j|\omega_{P'_4}|)\subseteq
\bigcup_{j=0}^{j_2}(\omega_{P_4}+j|\omega_{P_4}|).$$
We will also use the ``easy'' variant of this size, defined by the formula

$$\size_{3,e, \vec{\P}_{3,4}}(f_3):=
\sup_{\vec{P}\in \vec{\P}_{3,4}}
\frac{1}{|I_{\vec{P}}|}
\int_{\R}
|f_3(x)| \widetilde{\chi}_{I_{\vec{P}}}^C(x)
\chi_{\{ x : N_1(x)\in \omega_{P_3} \} }
\chi_{\{ x : N_2(x)\in \omega_{P_4} \} } dx.
$$

We also define the energy of $f_3$ by

\begin{equation}
\energy_{3,\vec{\P}_{3,4}}(f_3):=
\sum_{m_1=0}^{j_1}
\sum_{m_2=0}^{j_2}
\energy^{m_1, m_2}_{3,\vec{\P}_{3,4}}(f_3)
\end{equation}
while $\energy^{m_1, m_2}_{3,\vec{\P}_{3,4}}(f_3)$ is defined by

$$\energy^{m_1, m_2}_{3,\vec{\P}_{3,4}}(f_3):=
\sup_{n\in\Z}\sup_{\vec{\bf{D}}}
\frac{2^n}{M^2}
\left(\sum_{\vec{P'}\in\vec{\bf{D}}} |I_{\vec{P'}}|\right)
$$
where $\vec{\bf{D}}$ ranges over all collections of vector-tiles of dimension
$2$ $\vec{P'}\in \vec{\overline{\P}}_{3,4}$ for which there exists
$\vec{P}\in\vec{\P}_{3,4}$ with 
$\vec{P'}\in \vec{\overline{\P}}_{3,4}(\vec{P})$, having the property
that their corresponding parallelepipeds 
$I_{\vec{P'}}\times (\omega_{P'_3}+m_1|\omega_{P'_3}|)\times
(\omega_{P'_4}+m_2|\omega_{P'_4}|)$ are all disjoint in $\R^3$ and such that

$$\int_{\R}|f_3(x)|\widetilde{\chi}_{I_{\vec{P'}}}^C(x)
\chi_{\{ x : N_1(x)\in (\omega_{P'_3}+m_1|\omega_{P'_3}|) \}}
\cdot
\chi_{\{ x : N_2(x)\in (\omega_{P'_4}+m_2|\omega_{P'_4}|) \}}
dx \geq \frac{2^n}{M^2}|I_{\vec{P'}}|
$$
where $C>0$ in the above is a big constant which may vary from time to time 
depending on which $L^p$ estimates we are proving.

\end{definition}

The following general Proposition will play an important role in our estimates.

\begin{proposition}\label{use}
Let $\vec{\P}\subseteq \vec{\overline{\P}}$ be a finite collection of 
vector-tiles of dimension $4$ and let $f_1, f_2, f_3$ be fixed functions.
Then,

$$\left|\Lambda_{\vec{\P}}(f_1, f_2, f_3)\right|\lesssim
$$

\begin{equation}\label{19}
\left(\size_{1,\P_1}(f_1)\right)^{a}
\left(\size_{2,\P_2}(f_2)\right)^{a}
\left(\size_{3,\vec{\P}_{3,4}}(f_3)\right)^{b}\cdot
\end{equation}

$$\left(\energy_{1,\P_1}(f_1)\right)^{1-a}
\left(\energy_{2,\P_2}(f_2)\right)^{1-a}
\left(\energy_{3,\vec{\P}_{3,4}}(f_3)\right)^{1-b}$$
for any $0<a, b < 1$ with $a + 2 b = 1$, with the implicit
constants depending on $a, b, M$. Moreover, if for any
$\vec{P}, \vec{P'}\in \vec{\P}_{3,4}$ with 
$|I_{\vec{P}}|\neq|I_{\vec{P'}}| $one has 
$I_{\vec{P}}\cap I_{\vec{P'}} = \emptyset$, then the inequality (\ref{19})
holds even if one replaces $\size_{3,\vec{\P}_{3,4}}(f_3)$ with the smaller
quantity $\size_{3, e, \vec{\P}_{3,4}}(f_3)$.
\end{proposition}

\section{Upper bounds for sizes and energies}

The proof of Proposition \ref{use} will be postponed for a while. 
Until then, we will take advantage of it.
Clearly, in order to be able to use this Proposition \ref{use} effectively,
we need some estimates on sizes and energies.
The following two Lemmas have been proven in \cite{cct}, 
\cite{mtt:fourierbiest}.

\begin{lemma}\label{e12b}
Let $\vec{\P}$ be an arbitrary collection of vector-tiles of dimension $4$
and $f_i\in L^2(\R)$ for $i=1,2$. Then, we have

\begin{equation}
\energy_{i,\P_i}(f_i)\lesssim \|f_i\|_2
\end{equation}
for $i=1,2$.
\end{lemma}

\begin{lemma}\label{s12b}
Let $i=1,2$, $E_i$ be a set of finite measure, $f_i$ be a function in $X(E_i)$
and let also $\vec{\P}$ as before. Then, we have

\begin{equation}
\size_{i,\P_i}(f_i)\lesssim\sup_{\vec{P}\in\vec{\P}}
\frac{\int_{E_i}\widetilde{\chi}_{I_{\vec{P}}}^C(x) dx}{|I_{\vec{P}}|}
\end{equation}
for $i=1,2$ and any $C>0$ with the implicit constant depending on $C$.
\end{lemma}

The following Lemma follows immediately from Definition \ref{se34}.

\begin{lemma}\label{s3b}
Let $\vec{\P}_{3,4}\subseteq \vec{\overline{\P}}_{3,4}$ be as in 
Definition \ref{se34}, $E_3$ be a set with finite measure and $f_3\in X(E_3)$.
Then,

\begin{equation}
\size_{3,\vec{\P}_{3,4}}(f_3)\lesssim\sup_{\vec{P}\in\vec{\P}_{3,4}}
\sup_{\vec{P'}\in\vec{\overline{\P}}_{3,4}(\vec{P})}
\frac{\int_{E_3}\widetilde{\chi}_{I_{\vec{P'}}}^C(x) dx}{|I_{\vec{P'}}|}
\end{equation}
and similarly,

\begin{equation}
\size_{3,e, \vec{\P}_{3,4}}(f_3)\lesssim\sup_{\vec{P}\in\vec{\P}_{3,4}}
\frac{\int_{E_3}\widetilde{\chi}_{I_{\vec{P}}}^C(x) dx}{|I_{\vec{P}}|}
\end{equation}
where $C>0$ is the same constant which appeared in Definition \ref{se34}.
\end{lemma}

Finally, we also have

\begin{lemma}\label{e3b}
Let $\vec{\P}_{3,4}\subseteq \vec{\overline{\P}}_{3,4}$ be as before and
$f_3$ be an $L^1(\R)$ function. Then, one has

\begin{equation}\label{20}
\energy_{3,\vec{\P}_{3,4}}(f_3)\lesssim \|f_3\|_1.
\end{equation}
\end{lemma}

\begin{proof} Clearly,  from Definition \ref{se34} it is enough to show that
for any $0\leq m_1\leq j_1$ and $0\leq m_2\leq j_2$ the energies 
$\energy^{m_1, m_2}_{3,\vec{\P}_{3,4}}$ satisfy the estimate (\ref{20}).

Fix now $n$ and $\vec{\bf{D}}$ so that the suppremum is attained in the
definition of $\energy^{m_1, m_2}_{3,\vec{\P}_{3,4}}$. Then, we can write

$$\energy^{m_1, m_2}_{3,\vec{\P}_{3,4}}(f_3)\sim
\frac{2^n}{M^2}
\left(\sum_{\vec{P'}\in\vec{\bf{D}}} |I_{\vec{P'}}|\right)
$$

$$\lesssim\sum_{\vec{P'}\in\vec{\bf{D}}}
\int_{\R}|f_3(x)|\widetilde{\chi}_{I_{\vec{P'}}}^C(x)
\chi_{\{ x : N_1(x)\in (\omega_{P'_3}+m_1|\omega_{P'_3}|) \}}
\cdot
\chi_{\{ x : N_2(x)\in (\omega_{P'_4}+m_2|\omega_{P'_4}|) \}}
dx.
$$
Since we also know that the corresponding 
parallelepipeds 
$I_{\vec{P'}}\times (\omega_{P'_3}+m_1|\omega_{P'_3}|)\times
(\omega_{P'_4}+m_2|\omega_{P'_4}|)$ are all disjoint in $\R^3$, the claim 
follows from an argument similar to the one used in Proposition 3.1 in 
\cite{laceyt3}.

\end{proof}

\section{Proof of Theorem \ref{discrete} near the vertex $A_1$}

Fix $\vec{\P}$ a finite collection of vector-tiles of dimension $4$ as usual.
 We will show that the trilinear form $\Lambda_{\vec{\P}}$ is of restricted
type $\alpha$ for admissible $3$-tuples $(\alpha_1, \alpha_2, \alpha_3)$
arbitrarily close to $A_1$, so that their bad index is $3$. 

Fix $\alpha$ as above with $\alpha_1=\alpha_2$ and let $E_1, E_2, E_3$
be sets of finite measure. By scaling invariance, we can assume that $|E_3|=1$.
We therefore need to find a major set $E'_3\subseteq E_3$ so that

\begin{equation}\label{21}
\left|\Lambda_{\vec{\P}}(f_1, f_2, f_3)\right|\lesssim
|E_1|^{\alpha_1}
|E_2|^{\alpha_2}
|E_3|^{\alpha_3}
\end{equation}
for all functions $f_1\in X(E_1)$, $f_2\in X(E_2)$, $f_3\in X(E'_3)$.

Define the exceptional set $\Omega$ by

$$\Omega:= \bigcup_{j=1}^2
\left\{M(\frac{\chi_{E_j}}{|E_j|}) > C\right\}
$$
where $M$ is the classical Hardy-Littlewood maximal operator \cite{stein}.
Clearly, $|\Omega|<1/2$ if $C$ is a sufficiently large constant. Then if we set
$E'_3:=E_3\setminus\Omega$, $E'_3$ is a major subset of $E_3$. 

Let then $f_j\in X(E'_j)$ for $j=1,2,3$ and decompose the set $\vec{\P}$
as 

$$\vec{\P}= \bigcup_{d\geq 0}\vec{\P}_d$$
where $\vec{\P}_d$ is the set of all $\vec{P}\in\vec{\P}$ having the property
that 

$$1+\frac{\dist(I_{\vec{P}},\Omega^c)}{|I_{\vec{P}}|}\sim 2^d.$$
From the definition of $\Omega$ we have

$$\frac{1}{|I_{\vec{P}}|}
\int_{E_j}\widetilde{\chi}_{I_{\vec{P}}}\lesssim 2^d |E_j|
$$
for $j=1,2$ whenever $\vec{P}\in\vec{\P}_d$ and also

$$\frac{1}{|I_{\vec{P}}|}
\int_{E'_3}\widetilde{\chi}_{I_{\vec{P}}}\lesssim 2^{-Cd}
$$
whenever $\vec{P}\in\vec{\P}_d$ where $C>$ is an arbitrarily big constant.

If $0\leq d\leq 5$, it follows from Lemmas \ref{e12b} $-$ \ref{e3b} that

$$\size_{1, \P_{d,1}}(f_1)\lesssim |E_1|$$

$$\size_{2, \P_{d,2}}(f_2)\lesssim |E_2|$$

$$\size_{3, \vec{\P}_{d,3,4}}(f_3)\lesssim 1$$
and also that

$$\energy_{1, \P_{d,1}}(f_1)\lesssim |E_1|^{1/2}$$

$$\energy_{2, \P_{d,2}}(f_2)\lesssim |E_2|^{1/2}$$

$$\energy_{3, \vec{\P}_{d,3,4}}(f_3)\lesssim 1.$$
Similarly, for $d\geq 5$ using the same Lemmas we obtain

$$\size_{1, \P_{d,1}}(f_1)\lesssim 2^d|E_1|$$

$$\size_{2, \P_{d,2}}(f_2)\lesssim 2^d|E_2|$$

$$\size_{3, \vec{\P}_{d,3,4}}(f_3)\lesssim 2^{-Cd}$$
for any big number $C>0$, while the energies satisfy the same bounds as before.
Since the vector-tiles in $\vec{\P}_d$ (for $d\geq 5$) have the property 
that $\vec{P}, \vec{P'}\in\vec{\P}$ with 
$|I_{\vec{P}}|\neq |I_{\vec{P'}}|$ implies 
$I_{\vec{P}}\cap I_{\vec{P'}}=\emptyset$, it follows from 
Proposition \ref{use} that

$$\left|\Lambda_{\vec{\P}}(f_1, f_2, f_3)\right|
\lesssim |E_1|^a |E_2|^a |E_1|^{\frac{1-a}{2}}|E_2|^{\frac{1-a}{2}}+
$$

$$
\sum_{d\geq 5}(2^d|E_1|)^a(2^d|E_2|)^a 2^{-Cdb}
|E_1|^{\frac{1-a}{2}}|E_2|^{\frac{1-a}{2}}=
|E_1|^{\frac{1+a}{2}}|E_2|^{\frac{1+a}{2}}
\left(1+\sum_{d\geq 5}2^{2ad-Cdb}\right),$$
where $a$ and $b$ are as in Proposition \ref{use}.

Now if we pick $a$ close to $1$ so that $1/2+a/2=\alpha_1$, we have to define
$b$ by $b:=(1-a)/2$ and if we then choose $C$ big enough the expression above
becomes $|E_1|^{\alpha_1}|E_2|^{\alpha_2}$ which is the desired (\ref{21}).

\section{Proof of Theorem \ref{discrete} near the vertices $A_2$ and $A_3$}

Fix $\vec{\P}$ a finite collection of vector-tiles of dimension $4$ 
and let $\alpha$ an admissible $3$-tuple arbitrarily close to $A_2$ with
$\alpha_3=1$ so that the bad index is $1$ (the case of $A_3$ is similar,
by the symmetry of the form).

Let now $E_1, E_2, E_3$
be sets of finite measure. By scaling invariance, we can assume that $|E_1|=1$.
We then need to find a major subset $E'_1\subseteq E_1$ so that

\begin{equation}
\left|\Lambda_{\vec{\P}}(f_1, f_2, f_3)\right|\lesssim
|E_1|^{\alpha_1}
|E_2|^{\alpha_2}
|E_3|^{\alpha_3}
\end{equation}
for all functions $f_1\in X(E'_1)$, $f_2\in X(E_2)$, $f_3\in X(E_3)$.

Define as before the exceptional set $\Omega$ by

$$\Omega:= \bigcup_{j=2}^3
\left\{M(\frac{\chi_{E_j}}{|E_j|}) > C\right\}
$$
where $M$ is the Hardy-Littlewood maximal operator
and observe that $|\Omega|<1/2$ if $C$ is a sufficiently large constant. 
Then if we set
$E'_1:=E_1\setminus\Omega$, $E'_1$ is a major subset of $E_1$. 

Decompose again the set $\vec{\P}$
as

$$\vec{\P}= \bigcup_{d\geq 0}\vec{\P}_d.$$

From the definition of $\Omega$ we have

$$\frac{1}{|I_{\vec{P}}|}
\int_{E_j}\widetilde{\chi}_{I_{\vec{P}}}\lesssim 2^d |E_j|
$$
for $j=2,3$ whenever $\vec{P}\in\vec{\P}_d$ and also

$$\frac{1}{|I_{\vec{P}}|}
\int_{E'_1}\widetilde{\chi}_{I_{\vec{P}}}\lesssim 2^{-Cd}
$$
whenever $\vec{P}\in\vec{\P}_d$ where $C>$ is an arbitrarily big constant.

As before if $0\leq d\leq 5$, it follows from Lemmas \ref{e12b} $-$ \ref{e3b} 
that

$$\size_{1, \P_{d,1}}(f_1)\lesssim 1$$

$$\size_{2, \P_{d,2}}(f_2)\lesssim |E_2|$$

$$\size_{3, \vec{\P}_{d,3,4}}(f_3)\lesssim |E_3|$$
and also that

$$\energy_{1, \P_{d,1}}(f_1)\lesssim 1$$

$$\energy_{2, \P_{d,2}}(f_2)\lesssim |E_2|^{1/2}$$

$$\energy_{3, \vec{\P}_{d,3,4}}(f_3)\lesssim |E_3|.$$
Similarly, for $d\geq 5$ using the same Lemmas we obtain

$$\size_{1, \P_{d,1}}(f_1)\lesssim 2^{-Cd}$$

$$\size_{2, \P_{d,2}}(f_2)\lesssim 2^d|E_2|$$

$$\size_{3, \vec{\P}_{d,3,4}}(f_3)\lesssim 2^{d}|E_3|$$
for any big number $C>0$, while the energies satisfy the same bounds as before.
It follows again from 
Proposition \ref{use} that

$$\left|\Lambda_{\vec{\P}}(f_1, f_2, f_3)\right|
\lesssim |E_2|^a |E_3|^b |E_2|^{\frac{1-a}{2}}|E_3|^{1-b}+
$$

$$
\sum_{d\geq 5}
2^{-Cda}
(2^d|E_2|)^a (2^d|E_3|)^b 
|E_2|^{\frac{1-a}{2}} |E_3|^{1-b}=
|E_2|^{\frac{1+a}{2}}|E_3|
\left(1+\sum_{d\geq 5}2^{ad+bd-Cad}\right),$$
where $a$ and $b$ are as in Proposition \ref{use}.

Now if we pick $a$ close to $1$ so that $1/2+a/2=\alpha_2$, we have to define
$b$ by $b:=(1-a)/2$ and if we then choose $C$ big enough the expression above
becomes $|E_2|^{\alpha_2}|E_3|$ and this completes the proof.

\section{Phase space decompositions}

In order to complete our estimates, it remains to prove Proposition \ref{use}.
This will be accomplished with the help of certain combinatorial Lemmas. 
The first one is standard and it appeared in \cite{mtt:bicarleson}.

\begin{lemma}\label{dec12l}
Let $j=1,2$, $n\in\Z$, $\P'_j\subseteq \P_j$ be collection of tiles,
$f_1, f_2$ be two functions and suppose that

$$\size_{j,\P'_j}(f_j)\leq 2^{-n}\energy_{j,\P_j}(f_j).$$
Then, we can decompose $\P'_j$ as $\P'_j = \P''_j\cup \P'''_j$ such that

$$\size_{j,\P''_j}(f_j)\leq 2^{-n-1}\energy_{j,\P_j}(f_j)$$
and $\P'''_j$ can be written as a disjoint union of trees in $\T_j$ such that

$$\sum_{T\in\T_j} |I_T| \lesssim 2^{2n}.$$
\end{lemma}

By iterating this Lemma \ref{dec12l} one obtains 
(see again \cite{mtt:bicarleson}).

\begin{corollary}\label{dec12c}
With the same notations as in Lemma \ref{dec12l}, there exists a partition

$$\P_j = \bigcup_{n\in\Z}\P_j^n$$
where for each $n\in\Z$ one has

$$\size_{j,\P_j^n}(f_j)\leq 
\min (2^{-n}\energy_{j,\P_j}(f_j), \size_{j,\P_j}(f_j)).$$
Also, we can cover $\P_j^n$ by a collection $\T^n_j$ of trees such that

$$\sum_{T\in\T_j^n} |I_T| \lesssim 2^{2n}.$$
\end{corollary}

We will also need

\begin{lemma}\label{dec3l}
Let $\vec{\P'}_{3,4}\subseteq \vec{\P}_{3,4}$ be collections of vector-tiles
of dimension $2$, $f_3$ be a function and suppose that

\begin{equation}\label{22}
\size_{3,\vec{\P'}_{3,4}}(f_3)\leq 2^{-n}\energy_{3,\vec{\P}_{3,4}}(f_3).
\end{equation}
Then, we can decompose $\vec{\P'}_{3,4}$ as 
$\vec{\P'}_{3,4} = \vec{\P''}_{3,4}\cup \vec{\P'''}_{3,4}$ such that

\begin{equation}\label{23}
\size_{3,\vec{\P''}_{3,4}}(f_3)\leq 2^{-n-1}\energy_{3,\vec{\P}_{3,4}}(f_3)
\end{equation}
and $\vec{\P'''}_{3,4}$ can be written as a disjoint union of trees
$\vec{\T}_{3,4}$ such that

\begin{equation}\label{24}
\sum_{\vec{T}\in \vec{\T}_{3,4}} |I_{T_3}\cap I_{T_4}|\lesssim 2^n.
\end{equation}
Moreover, if for any $\vec{P}, \vec{P'}\in \vec{\P}_{3,4}$ with
$|I_{\vec{P}}|\neq |I_{\vec{P'}}|$ one has 
$I_{\vec{P}}\cap I_{\vec{P'}} = \emptyset$, then the above statement holds 
even if one replaces ``$\size_3...$'' with the smaller quantity 
``$\size_{3, e}...$''.
\end{lemma}

\begin{proof} First, we consider all the one vector-tile collections
$\{\vec{P}\}$ with $\vec{P}\in \vec{\P'}_{3,4}$, having the property
that

\begin{equation}\label{25}
\size_{3, \{\vec{P}\}}(f_3) >
2^{-n-1}\energy_{3,\vec{\P}_{3,4}}(f_3).
\end{equation}
From Definition \ref{se34} it follows that for each such a $\vec{P}$, 
there exists 
$\vec{P'}\in\vec{\overline{\P}}_{3,4}(\vec{P})$ such that

\begin{equation}\label{26}
\frac{1}{|I_{\vec{P'}}|}
\int_{\R}
|f_3(x)| \widetilde{\chi}_{I_{\vec{P'}}}^C(x)
\chi_{\{ x : N_1(x)\in \cup_{j=0}^{j_1}(\omega_{P'_3}+j|\omega_{P'_3}|) \}}
\cdot
\chi_{\{ x : N_2(x)\in \cup_{j=0}^{j_2}(\omega_{P'_4}+j|\omega_{P'_4}|) \}}
dx > 
\end{equation}

$$2^{-n-1}\energy_{3,\vec{\P}_{3,4}}(f_3).$$
Clearly, for any such a $\vec{P'}$ there exist two
indices $0\leq m_1(\vec{P'})\leq j_1$ and $0\leq m_2(\vec{P'})\leq j_2$
so that
\begin{equation}\label{27}
\frac{1}{|I_{\vec{P'}}|}
\int_{\R}
|f_3(x)| \widetilde{\chi}_{I_{\vec{P'}}}^C(x)
\chi_{\{ x : N_1(x)\in (\omega_{P'_3}+m_1(\vec{P'}) |\omega_{P'_3}|) \}}
\cdot
\chi_{\{ x : N_2(x)\in (\omega_{P'_4}+m_2(\vec{P'}) |\omega_{P'_4}|) \}}
dx > 
\end{equation}

$$\frac{2^{-n-1}}{M^2}\energy_{3,\vec{\P}_{3,4}}(f_3).$$
Now for any fixed indices $0\leq m_1\leq j_1$ and 
$0\leq m_2\leq j_2$ we denote by $\vec{\overline{\P}}_{3,4}(m_1, m_2)$
the set of all $\vec{P'}$ as above with the property that 
$m_1(\vec{P'})=m_1$ and $m_2(\vec{P'})=m_2$. We then introduce an ordering
$<^c_{m_1, m_2}$ on $\vec{\overline{\P}}_{3,4}(m_1, m_2)$ and write
$\vec{P'}<^c_{m_1, m_2}\vec{P''}$ if and only if

$$I_{\vec{P'}}\times (\omega_{P'_3}+m_1|\omega_{P'_3}|) <^c
I_{\vec{P''}}\times (\omega_{P''_3}+m_1|\omega_{P''_3}|)$$
and

$$I_{\vec{P'}}\times (\omega_{P'_4}+m_1|\omega_{P'_4}|) <^c
I_{\vec{P''}}\times (\omega_{P''_4}+m_1|\omega_{P''_4}|)$$
where $<^c$ is C.Fefferman's classical ordering between tiles defined in
Definition \ref{order}.

Denote by $\vec{\overline{\P}}_{3,4}^{\text{max}}(m_1, m_2)$
the set of all $\vec{P'}\in \vec{\overline{\P}}_{3,4}(m_1, m_2)$
which are maximal with respect to this ordering $<^c_{m_1, m_2}$ defined
before.

Fix a generic $\vec{P'}\in \vec{\overline{\P}}_{3,4}^{\text{max}}(m_1, m_2)$
and consider now all the vector-tiles $\vec{P}\in\vec{\P'}_{3,4}$ 
with the property
that

\begin{equation}\label{28}
I_{\vec{P}}\subsetneq I_{\vec{P'}},
\end{equation}

\begin{equation}\label{29}
\bigcup_{j=0}^{j_1}(\omega_{P'_3}+j|\omega_{P'_3}|)\subseteq
\bigcup_{j=0}^{j_1}(\omega_{P_3}+j|\omega_{P_3}|)
\end{equation}
and

\begin{equation}\label{30}
\bigcup_{j=0}^{j_2}(\omega_{P'_4}+j|\omega_{P'_4}|)\subseteq
\bigcup_{j=0}^{j_2}(\omega_{P_4}+j|\omega_{P_4}|).
\end{equation}
This set can be reorganized as a union of subsets denoted by
$\vec{T}_{\vec{P'}}^{s_1, s_2}$ for $0\leq s_1\leq j_1$ and
$0\leq s_2\leq j_2$ where $\vec{T}_{\vec{P'}}^{s_1, s_2}$ contains all the
vector-tiles $\vec{P}$ satisfying (\ref{28}), (\ref{29}), (\ref{30})
which have also the additional property that

$$\omega_{P'_3}+m_1|\omega_{P'_3}|\subseteq \omega_{P_3}+s_1|\omega_{P_3}|$$
and

$$\omega_{P'_4}+m_2|\omega_{P'_4}|\subseteq \omega_{P_4}+s_2|\omega_{P_4}|.$$
It is easy to see that all this collections 
$\vec{T}_{\vec{P'}}^{s_1, s_2}$ are vector-trees $\vec{T}$ of dimension $2$,
for which

\begin{equation}\label{31}
I_{T_3}= I_{T_4}=I_{\vec{P'}}.
\end{equation}
Also, by construction, it follows that all the $\vec{P}$'s satisfying
(\ref{25}) have been selected this way. We collect now all these
vector-trees into a set named $\vec{\P'''}_{3,4}$ and it is easy to observe
that if we set 
$\vec{\P''}_{3,4}:= \vec{\P'}_{3,4}\setminus \vec{\P'''}_{3,4}$, then
(\ref{23}) is satisfied. 

It remains to prove (\ref{24}) only. From (\ref{31}), it is clearly enough
to demonstrate that

\begin{equation}\label{32}
\sum_{m_1=0}^{j_1}
\sum_{m_2=0}^{j_2}
\sum_{\vec{\P'}\in 
\vec{\overline{\P}}_{3,4}^{\text{max}}(m_1, m_2)}
|I_{\vec{P'}}|\lesssim 2^n.
\end{equation}
By using Definition \ref{se34} and our selection algorithm, we see that
the left hand side of (\ref{32}) can be majorized by

$$2^n (\energy_{3,\vec{\P}_{3,4}}(f_3))^{-1}
\sum_{m_1=0}^{j_1}
\sum_{m_2=0}^{j_2}
\energy^{m_1, m_2}_{3,\vec{\P}_{3,4}}(f_3)\lesssim 2^n
$$
and this completes the proof.

\end{proof}

As before, by iterating this Lemma \ref{dec3l}, one obtains

\begin{corollary}\label{dec3c}
With the same notations as in Lemma \ref{dec3l}, there exists a partition

$$\vec{\P}_{3,4} = \bigcup_{n\in\Z}\vec{\P}^n_{3,4}$$
where for each $n\in\Z$ one has

$$\size_{3,\vec{\P}^n_{3,4}}\leq 
\min (2^{-n} \energy_{3,\vec{\P}_{3,4}}, \size_{3,\vec{\P}_{3,4}}).
$$
Also, we can cover $\vec{\P}^n_{3,4}$ by a collection $\vec{\T}^n_{3,4}$
of trees such that

$$\sum_{\vec{T}\in\vec{\T}^n_{3,4}}
|I_{T_3}\cap I_{T_4}|\lesssim 2^n.$$
\end{corollary}

Fix now a finite collection of vector-tiles of dimension 
$4$ $\vec{\P}\subseteq \vec{\overline{\P}}$ and $f_1, f_2, f_3$ three
functions. Consider also its projected collections $\P_1, \P_2$ and
$\vec{\P}_{3,4}$. Pick now $\vec{P}\in\vec{\P}$ arbitrary. By using
Corollaries \ref{dec12c} and \ref{dec3c}, there exist $n_1, n_2, n_3\in\Z$
and trees $T_1\in\T_1^{n_1}$, $T_2\in\T_2^{n_2}$ and 
$\vec{T}_{3,4}\in\vec{\T}^{n_3}_{3,4}$ so that $P_1\in T_1$, $P_2\in\T_2$
and $(P_3, P_4)\in\vec{T}_{3,4}$. Then, consider all the other vector-tiles
$\vec{P'}\in \vec{\P}$ having the same property that
$P'_1\in T_1$, $P'_2\in T_2$ and $(P'_3, P'_4)\in\vec{T}_{3,4}$.
Clearly, they form a vector-tree of dimension $4$ and as a consequence,
our initial collection $\vec{\P}$ can be written as a disjoint union
of such trees $\vec{T}$. To estimate our trilinear form, the standard way 
is to first understand the contribution of a single tree. This is the scope 
of the next Lemma.

\begin{lemma}\label{treeestimate}
Let $\vec{T}$ be one of the vector-trees constructed before. Then,

\begin{equation}\label{te}
\left|\Lambda_{\vec{T}}(f_1, f_2, f_3)\right|\lesssim
\size_{1, T_1}(f_1)
\size_{2, T_2}(f_2)
\size_{3, \vec{T}_{3,4}}(f_3)
|I_{T_1}\cap I_{T_2}\cap I_{T_3}\cap I_{T_4}|.
\end{equation}
\end{lemma}
The proof of this important Lemma wil be presented later on. In the meantime 
we will take it for granted in order to complete the proof of 
Proposition \ref{use}.

\section{Proof of Proposition \ref{use}}

We denote, for simplicity, by $S_1, S_2, S_3$ and $E_1, E_2, E_3$ the three
sizes and energies which appear in the inequality (\ref{19}). Then, by using
the Corollaries \ref{dec12c} and \ref{dec3c} and also Lemma 11.5, we can write

$$\left|\Lambda_{\vec{P}}(f_1, f_2, f_3)\right|
= E_1 E_2 E_3 \left|\Lambda_{\vec{P}}
(\frac{f_1}{E_1}, \frac{f_2}{E_2}, \frac{f_3}{E_3})\right|\lesssim
$$

\begin{equation}\label{33}
\sum_{n_1, n_2, n_3} 2^{-n_1} 2^{-n_2} 2^{-n_3}
\sum_{\vec{T}\in\vec{\T}_{n_1, n_2, n_3}}
|I_{T_1}\cap I_{T_2}\cap I_{T_3}\cap I_{T_4}|,
\end{equation}
where $\vec{\T}_{n_1, n_2, n_3}$ is the set of all vector-trees $\vec{T}$
described at the end of the previous section and the summation goes over
the indices $n_1, n_2, n_3 \in \Z$ having the property that

\begin{equation}\label{constraints}
2^{-n_j}\lesssim \frac{S_j}{E_j}
\end{equation}
for $j=1,2,3$. Let us now recall the fact that there are actually two degrees
of freedom in our vector-tiles of dimension $4$ 
$\vec{P} = (P_1, P_2, P_3, P_4)$, since once we fix the tiles $P_1$ and $P_2$
then $P_3$ and $P_4$ are uniquely determined and viceversa. As a consequence
of this fact, it is not difficult to observe that

$$\sum_{\vec{T}\in\vec{\T}_{n_1, n_2, n_3}}
|I_{T_1}\cap I_{T_2}\cap I_{T_3}\cap I_{T_4}|
\lesssim \left(\sum_{T\in \T^{n_1}_1} |I_T|\right)
\left(\sum_{T\in \T^{n_2}_2} |I_T|\right)$$
and also that

$$\sum_{\vec{T}\in\vec{\T}_{n_1, n_2, n_3}}
|I_{T_1}\cap I_{T_2}\cap I_{T_3}\cap I_{T_4}|
\lesssim \sum_{\vec{T}\in\vec{\T}^{n_3}_{3,4}}|I_{T_3}\cap I_{T_4}|.
$$
In particular, using the same Corollaries, it follows that

\begin{equation}
\sum_{\vec{T}\in\vec{\T}_{n_1, n_2, n_3}}
|I_{T_1}\cap I_{T_2}\cap I_{T_3}\cap I_{T_4}|
\lesssim \min (2^{2 n_1} 2^{ 2 n_2},\, 2^{n_3})
\end{equation}
and this implies that

\begin{equation}\label{34}
\sum_{\vec{T}\in\vec{\T}_{n_1, n_2, n_3}}
|I_{T_1}\cap I_{T_2}\cap I_{T_3}\cap I_{T_4}|
\lesssim 2^{2 n_1 a'} 2^{2 n_2 a'} 2^{n_3 b'}
\end{equation}
for any $a', b'\in [0,1]$ with $a'+ b' = 1$.

Using (\ref{34}) in (\ref{33}) we estimate our trilinear form by

$$E_1 E_2 E_3
\sum_{n_1, n_2, n_3}
2^{-n_1(1-2a')}
2^{-n_2(1-2a')}
2^{-n_3( 1-b')}\lesssim
$$

$$E_1 E_2 E_3
\left(\frac{S_1}{E_1}\right)^{1-2a'}
\left(\frac{S_2}{E_2}\right)^{1-2a'}
\left(\frac{S_3}{E_3}\right)^{1-b'} =
S_1^{1-2a'} S_2^{1-2a'} S_3^{1-b'} E_1^{2a'} E_2^{2a'} E_3^{b'},$$
if we assume that $a' < 1/2$ and $b' < 1$ and use the constraints
(\ref{constraints}).

Now we just have to define $a:= 1-2a'$ and $b:= 1-b'$ and to observe that
$a+2b=1$, in order to complete the proof.

\section{Proof of Lemma \ref{treeestimate}}

We are left with proving Lemma \ref{treeestimate}. Fix $\vec{T}$ a 
vector-tree of dimension $4$ as there and denote by
$I_{\vec{T}}: I_{T_1}\cap I_{T_2}\cap I_{T_3}\cap I_{T_4}$. Clearly,
any $\vec{P}\in\vec{T}$ has the property that 
$I_{\vec{P}}\subseteq I_{\vec{T}}$. 

Consider now $\J$ the collection
of all dyadic intervals $J\subseteq\R$ having the property that $3J$ does not
contain any $I_{\vec{P}}$ for $\vec{P}\in\vec{T}$ and that $J$ is maximal with
this property. Clearly, all this intervals $J\in\J$ are disjoint and 
their union is the whole real line $\R$. We also observe that if $J\in\J$
satisfies $J\cap 3I_{\vec{T}}\neq\emptyset$, then $J\subseteq 3I_{\vec{T}}$
and moreover, if $J\subseteq (3I_{\vec{T}})^c$ and $\vec{P}\in\vec{T}$,
then one has $|I_{\vec{P}}|\leq |J|$.
As a consequence of these, we can write

$$\Lambda_{\vec{T}}(f_1, f_2, f_3)=$$

$$\int_{\R}\left(
\sum_{\vec{P}\in \vec{T}}
\frac{1}{|I_{\vec{P}}|^{1/2}}
\langle f_1, \Phi_{P_1}\rangle
\langle f_2, \Phi_{P_2}\rangle
\Phi_{\vec{P}}
\chi_{\{x: N_1(x)\in\omega_{P_3}\}}
\chi_{\{x: N_2(x)\in\omega_{P_4}\}}
\right) f_3(x) dx=
$$

$$\sum_{J\subseteq 3I_{\vec{T}}}
\int_{J}\left(
\sum_{\vec{P}\in \vec{T}}
\frac{1}{|I_{\vec{P}}|^{1/2}}
\langle f_1, \Phi_{P_1}\rangle
\langle f_2, \Phi_{P_2}\rangle
\Phi_{\vec{P}}
\chi_{\{x: N_1(x)\in\omega_{P_3}\}}
\chi_{\{x: N_2(x)\in\omega_{P_4}\}}
\right) f_3(x) dx + 
$$

$$\sum_{J\subseteq (3I_{\vec{T}})^c}
\int_{J}\left(
\sum_{\vec{P}\in \vec{T}}
\frac{1}{|I_{\vec{P}}|^{1/2}}
\langle f_1, \Phi_{P_1}\rangle
\langle f_2, \Phi_{P_2}\rangle
\Phi_{\vec{P}}
\chi_{\{x: N_1(x)\in\omega_{P_3}\}}
\chi_{\{x: N_2(x)\in\omega_{P_4}\}}
\right) f_3(x) dx := 
$$

\begin{equation}\label{unudoi}
I + II.
\end{equation}
To be accurate, we should mention that the functions $\Phi_{\vec{P}}$
which appear in the above expressions, are in fact complex conjugates
of the previous functions $\Phi_{\vec{P}}$, considered in (\ref{discreteform}).

Term $II$ can be viewed as an error term and can be estimated by

$$\sum_{J\subseteq (3I_{\vec{T}})^c}
\sum_{\vec{P}\in \vec{T}}
\int_{J}
\frac{1}{|I_{\vec{P}}|^{1/2}}
\sum_{\vec{P}\in \vec{T}}
\frac{1}{|I_{\vec{P}}|^{1/2}}
|\langle f_1, \Phi_{P_1}\rangle|
|\langle f_2, \Phi_{P_2}\rangle|
|\Phi_{\vec{P}}|
\chi_{\{x: N_1(x)\in\omega_{P_3}\}}
\chi_{\{x: N_2(x)\in\omega_{P_4}\}}
|f_3(x)| dx\lesssim
$$

$$\size_{1, T_1}(f_1)
\size_{2, T_2}(f_2)
\sum_{J\subseteq (3I_{\vec{T}})^c}
\sum_{\vec{P}\in \vec{T}}
\int_{J}
\widetilde{\chi}_{I_{\vec{P}}}^m(x)
\chi_{\{x: N_1(x)\in\omega_{P_3}\}}
\chi_{\{x: N_2(x)\in\omega_{P_4}\}}
|f_3(x)| dx
$$
for any positive integer $m$. Then, the above expression can be further
estimated by

$$\size_{1, T_1}(f_1)
\size_{2, T_2}(f_2)
\sum_{J\subseteq (3I_{\vec{T}})^c}
\sum_{\vec{P}\in \vec{T}}
\size_{3, \vec{T}_{3,4}}(f_3)
|I_{\vec{P}}|
\left(\frac{\dist(I_{\vec{P}}, J)}{|I_{\vec{P}}|}\right)^{-m}=
$$

$$\size_{1, T_1}(f_1)
\size_{2, T_2}(f_2)
\size_{3, \vec{T}_{3,4}}(f_3)
\sum_{J\subseteq (3I_{\vec{T}})^c}
\sum_{k\in\Z: 2^k\leq |J|}
\sum_{\vec{P}: |I_{\vec{P}}|= 2^k}
|I_{\vec{P}}|
\left(\frac{\dist(I_{\vec{P}}, J)}{|I_{\vec{P}}|}\right)^{-m}=
$$

$$\size_{1, T_1}(f_1)
\size_{2, T_2}(f_2)
\size_{3, \vec{T}_{3,4}}(f_3)
\sum_{J\subseteq (3I_{\vec{T}})^c}
\sum_{k\in\Z: 2^k\leq |J|} 2^k
\sum_{\vec{P}: |I_{\vec{P}}|= 2^k}
\left(\frac{\dist(I_{\vec{P}}, J)}{|I_{\vec{P}}|}\right)^{-m}\lesssim
$$

$$\size_{1, T_1}(f_1)
\size_{2, T_2}(f_2)
\size_{3, \vec{T}_{3,4}}(f_3)
\sum_{J\subseteq (3I_{\vec{T}})^c}
\sum_{k\in\Z: 2^k\leq |J|} 2^k
\left(\frac{\dist(I_{\vec{T}}, J)}{|I_{\vec{T}}|}\right)^{-m}\lesssim
$$

$$\size_{1, T_1}(f_1)
\size_{2, T_2}(f_2)
\size_{3, \vec{T}_{3,4}}(f_3)
\sum_{J\subseteq (3I_{\vec{T}})^c}
|J|
\left(\frac{\dist(I_{\vec{T}}, J)}{|I_{\vec{T}}|}\right)^{-m}\lesssim
$$

\begin{equation}
\size_{1, T_1}(f_1)
\size_{2, T_2}(f_2)
\size_{3, \vec{T}_{3,4}}(f_3)
|I_{\vec{T}}|.
\end{equation}
It remains to estimate term $I$ in (\ref{unudoi}). 
At this moment, let us recall from the 
proof of Lemma \ref{dec3l} that we may assume that the projected
vector-tree of dimension $2$ $\vec{T}_{3,4}$, is ``of type $(s_1, s_2)$''
with top $\vec{P}_{\vec{T}} = (P_{\vec{T}, 3}, P_{\vec{T}, 4}) $,
meaning that the frequency intervals 
$(\omega_{P_3}+s_1|\omega_{P_3}|)_{\vec{P}\in\vec{T}}$ all contain
the top frequency $\omega_{P_{\vec{T}, 3}}$
and also that the intervals
$(\omega_{P_4}+s_2|\omega_{P_4}|)_{\vec{P}\in\vec{T}}$
all contain the top frequency $\omega_{P_{\vec{T}, 4}}$  . 
Recall also that $0\leq s_1\leq j_1$
 and $0\leq s_2\leq j_2$ and that $\max(j_1, j_2)=M$
and we will assume from this point on that $j_2 = M$ (the other situation being
of course, similar).

We have several cases.

\underline{Case $1$: $s_1\geq 2$ or $s_2\geq 2$.}

In this case, either the intervals $(\omega_{P_3})_{\vec{P}\in\vec{T}}$
or $(\omega_{P_4})_{\vec{P}\in\vec{T}}$ are all disjoint for different scales
and we can write

$$ I = \sum_{J\subseteq 3I_{\vec{T}}}
\sum_{|I_{\vec{P}}| > |J|}
\int_{J}
\frac{1}{|I_{\vec{P}}|^{1/2}}
\langle f_1, \Phi_{P_1}\rangle
\langle f_2, \Phi_{P_2}\rangle
\Phi_{\vec{P}}
\chi_{\{x: N_1(x)\in\omega_{P_3}\}}
\chi_{\{x: N_2(x)\in\omega_{P_4}\}}
f_3(x) dx + 
$$

$$\sum_{J\subseteq 3I_{\vec{T}}}
\sum_{|I_{\vec{P}}|\leq |J|}
\int_{J}
\frac{1}{|I_{\vec{P}}|^{1/2}}
\langle f_1, \Phi_{P_1}\rangle
\langle f_2, \Phi_{P_2}\rangle
\Phi_{\vec{P}}
\chi_{\{x: N_1(x)\in\omega_{P_3}\}}
\chi_{\{x: N_2(x)\in\omega_{P_4}\}}
f_3(x) dx:=
$$

$$I_1 + I_2.$$
Term $I_2$ can be estimated using the same argument we used when estimating
Term $II$, by

\begin{equation}
\size_{1, T_1}(f_1)
\size_{1, T_2}(f_2)
\size_{3, \vec{T}_{3,4}}(f_3)
\sum_{J\subseteq 3I_{\vec{T}}} |J|
\lesssim
\size_{1, T_1}(f_1)
\size_{1, T_2}(f_2)
\size_{3, \vec{T}_{3,4}}(f_3)|I_{\vec{T}}|.
\end{equation}

On the other hand, term $I_1$ is smaller than

\begin{equation}\label{35}
I_1\lesssim
\size_{1, T_1}(f_1)
\size_{1, T_2}(f_2)
\sum_{J\subseteq 3I_{\vec{T}}}
\int_J
\sum_{|I_{\vec{P}}| > |J|}
\widetilde{\chi}_{I_{\vec{P}}}^m(x)
\chi_{\{x: N_1(x)\in\omega_{P_3}\}}
\chi_{\{x: N_2(x)\in\omega_{P_4}\}}
|f_3(x)| dx
\end{equation}
for any positive integer $m$.

For every $J\in\J$, $J\subseteq 3I_{\vec{T}}$ consider the unique dyadic
interval $J\subseteq \tilde{J}$ so that $|\tilde{J}|= 2 |J|$. From the 
maximality of $J$, it follows that there exists $\vec{P^0}\in\vec{T}$
with $I_{\vec{P^0}}\subseteq 3\tilde{J}$.
Pick now a vector-tile $\vec{P}(\tilde{J})$ not necessarily in our tree
$\vec{T}$ such that $I_{\vec{P^0}}\subseteq I_{\vec{P}(\tilde{J})}$,
$|I_{\vec{P}(\tilde{J})} |= |\tilde{J}|$ and with

\begin{equation}
\omega_{P_{\vec{T}, 3}}\subseteq 
\omega_{P_3(\tilde{J})}+ s_1 |\omega_{P_3(\tilde{J})}|\subseteq
\omega_{P^0_3}+ s_1 |\omega_{P^0_3}|
\end{equation}
and

\begin{equation}
\omega_{P_{\vec{T}, 4}}\subseteq 
\omega_{P_4(\tilde{J})}+ s_2 |\omega_{P_4(\tilde{J})}|\subseteq
\omega_{P^0_4}+ s_2 |\omega_{P^0_4}|.
\end{equation}
In particular, this implies that

$$
\bigcup_{j=0}^{j_1}
(\omega_{P_3(\tilde{J})}+ j |\omega_{P_3(\tilde{J})}|)
\subseteq
\bigcup_{j=0}^{j_1}
(\omega_{P^0_3}+ j |\omega_{P^0_3}|)
$$
and

$$
\bigcup_{j=0}^{j_2}
(\omega_{P_4(\tilde{J})}+ j |\omega_{P_4(\tilde{J})}|)
\subseteq
\bigcup_{j=0}^{j_2}
(\omega_{P^0_4}+ j |\omega_{P^0_4}|)
$$
and also that

$$\bigcup_{j=0}^{j_1}
(\omega_{P_3} + j|\omega_{P_3}|)
\subseteq
\bigcup_{j=0}^{j_1}
(\omega_{P_3(\tilde{J})}+ j |\omega_{P_3(\tilde{J})}|)
$$
and

$$\bigcup_{j=0}^{j_2}
(\omega_{P_4} + j|\omega_{P_4}|)
\subseteq
\bigcup_{j=0}^{j_2}
(\omega_{P_4(\tilde{J})}+ j |\omega_{P_4(\tilde{J})}|)
$$
for any $\vec{P}\in\vec{T}$ with $|I_{\vec{P}}| > |J|$.
Since for different scales, the functions 
$\chi_{\{x: N_1(x)\in\omega_{P_3}\}}
\chi_{\{x: N_2(x)\in\omega_{P_4}\}}$
have disjoint supports, it follows that (\ref{35}) can be majorized by

$$\size_{1, T_1}(f_1)\size_{2, T_2}(f_2)\cdot$$

$$
\sum_{J\subseteq 3I_{\vec{T}}}
\int_{\R}
\widetilde{\chi}_{I_{\vec{P}(\tilde{J})}}^m(x)
\chi_{\{x: N_1(x)\in\bigcup_{j=0}^{j_1}
(\omega_{P_3(\tilde{J})}+ j |\omega_{P_3(\tilde{J})}|)\}}
\chi_{\{x: N_2(x)\in \bigcup_{j=0}^{j_2}
(\omega_{P_4(\tilde{J})}+ j |\omega_{P_4(\tilde{J})}|)\}}
|f_3(x)| dx\lesssim$$

$$\size_{1, T_1}(f_1)\size_{2, T_2}(f_2)
\size_{3, \vec{T}_{3,4}}(f_3)
\sum_{J\subseteq 3I_{\vec{T}}}
|I_{\vec{P}(\tilde{J})}|
\lesssim
\size_{1, T_1}(f_1)\size_{2, T_2}(f_2)
\size_{3, \vec{T}_{3,4}}(f_3)
|I_{\vec{T}}|,
$$
as desired.

\underline{Case $2$: $s_1\leq 1$ and $s_2\leq 1$.}

In this case, the intervals 
$(\omega_{P_3})_{\vec{P}\in\vec{T}}$
and $(\omega_{P_4})_{\vec{P}\in\vec{T}}$ are intersecting each other.
There are two subcases.

\underline{Subcase $2'$: $j_1\leq 5$.}

As in the previous case, Term $I$ can be estimated by 

$$I\lesssim I_1 + I_2$$
and it is enough to discuss $I_1$ only, $I_2$ being similar.
This time, we also know that the intervals
$(\omega_{P_2})_{\vec{P}\in\vec{T}}$ are disjoint for different scales,
while $(\omega_{P_1})_{\vec{P}\in\vec{T}}$ may intersect each other.
As a consequence, it is easy to see that the frequency intervals of the
functions $\Phi_{\vec{P}}$ given by $\supp\widehat{\Phi_{\vec{P}}}$
are all disjoint for different scales. Using this fact, an argument similar
to the one in \cite{laceyt3} and \cite{fefferman} allows us to estimate
Term $I_1$ by

$$
\sum_{J\subseteq 3I_{\vec{T}}}
(\sup_{J\subseteq I}
\frac{1}{|I|}
\int_I
|
\sum_{\vec{P}\in\vec{T}}\frac{1}{|I_{\vec{P}}|^{1/2}}
\langle f_1, \phi_{P_1}\rangle
\langle f_2, \phi_{P_2}\rangle
\phi_{\vec{P}}(x)| dx)\cdot
$$

$$\int_{\R}
\widetilde{\chi}_{I_{\vec{P}(\tilde{J})}}^m(x)
\chi_{\{x: N_1(x)\in\bigcup_{j=0}^{j_1}
(\omega_{P_3(\tilde{J})}+ j |\omega_{P_3(\tilde{J})}|)\}}
\chi_{\{x: N_2(x)\in \bigcup_{j=0}^{j_2}
(\omega_{P_4(\tilde{J})}+ j |\omega_{P_4(\tilde{J})}|)\}}
|f_3(x)| dx\lesssim$$

$$\sum_{J\subseteq 3I_{\vec{T}}}
(\sup_{J\subseteq I}
\frac{1}{|I|}
\int_I
|
\sum_{\vec{P}\in\vec{T}}\frac{1}{|I_{\vec{P}}|^{1/2}}
\langle f_1, \phi_{P_1}\rangle
\langle f_2, \phi_{P_2}\rangle
\phi_{\vec{P}}(x)| dx)
\size_{3,\vec{T}_{3,4}}(f_3) |I_{\vec{P}(\tilde{J})}|\lesssim
$$

$$\size_{3,\vec{T}_{3,4}}(f_3)
\sum_{J\subseteq 3I_{\vec{T}}}
(\sup_{J\subseteq I}
\frac{1}{|I|}
\int_I
|
\sum_{\vec{P}\in\vec{T}}\frac{1}{|I_{\vec{P}}|^{1/2}}
\langle f_1, \phi_{P_1}\rangle
\langle f_2, \phi_{P_2}\rangle
\phi_{\vec{P}}(x)| dx) |J|\lesssim
$$

$$\size_{3,\vec{T}_{3,4}}(f_3)
\|M(\sum_{\vec{P}\in\vec{T}}\frac{1}{|I_{\vec{P}}|^{1/2}}
\langle f_1, \phi_{P_1}\rangle
\langle f_2, \phi_{P_2}\rangle
\phi_{\vec{P}}(x) )\|_{L^1(3I_{\vec{T}})}\lesssim
$$

$$\size_{3,\vec{T}_{3,4}}(f_3)
\|M(\sum_{\vec{P}\in\vec{T}}\frac{1}{|I_{\vec{P}}|^{1/2}}
\langle f_1, \phi_{P_1}\rangle
\langle f_2, \phi_{P_2}\rangle
\phi_{\vec{P}}(x) )\|_{L^2(3I_{\vec{T}})}|I_{\vec{T}}|^{1/2}\lesssim
$$

$$\size_{3,\vec{T}_{3,4}}(f_3)
\|\sum_{\vec{P}\in\vec{T}}\frac{1}{|I_{\vec{P}}|^{1/2}}
\langle f_1, \phi_{P_1}\rangle
\langle f_2, \phi_{P_2}\rangle
\phi_{\vec{P}}(x)\|_{L^2}|I_{\vec{T}}|^{1/2}\lesssim
$$

\begin{equation}\label{36}
\size_{3,\vec{T}_{3,4}}(f_3)
(\sum_{\vec{P}\in\vec{T}}\frac{1}{|I_{\vec{P}}|^{1/2}}
|\langle f_1, \phi_{P_1}\rangle|
|\langle f_2, \phi_{P_2}\rangle|
|\langle h, \Phi_{\vec{P}}\rangle|)|I_{\vec{T}}|^{1/2}
\end{equation}
for some $h\in L^2(\R)$, $\|h\|_2 = 1$.
Then, (\ref{36}) can be further estimated by

$$\size_{3,\vec{T}_{3,4}}(f_3)
\size_{1, T_1}(f_1)
(\sum_{\vec{P}\in\vec{T}}
|\langle f_2, \phi_{P_2}\rangle|^2)^{1/2}
(\sum_{\vec{P}\in\vec{T}}
|\langle h, \Phi_{\vec{P}}\rangle|^2)^{1/2}|I_{\vec{T}}|^{1/2}\lesssim
$$

$$\size_{1, T_1}(f_1)\size_{1, T_2}(f_2)
\size_{3,\vec{T}_{3,4}}(f_3)|I_{\vec{T}}|$$
again by using the fact that the frequency intervals of the functions
$\Phi_{\vec{P}}$ are disjoint for different scales. Finally, we have

\underline{Subcase 2'': $j_1\geq 5$.}

This is the easiest situation since both collections of intervals
$(\omega_{P_1})_{\vec{P}\in\vec{T}}$ and
$(\omega_{P_2})_{\vec{P}\in\vec{T}}$ are disjoint for different scales.
Now it is not necessary to split Term $I$ as $I_1 + I_2$, but instead
simply write

$$I\lesssim
\sum_{\vec{P}\in\vec{T}}
|\langle f_1, \phi_{P_1}\rangle|
|\langle f_2, \phi_{P_2}\rangle|
\frac{1}{|I_{\vec{P}}|}
\int_{3 I_{\vec{T}}}
\widetilde{\chi}_{I_{\vec{P}}}^m(x)
\chi_{\{x: N_1(x)\in\omega_{P_3}\}}
\chi_{\{x: N_2(x)\in\omega_{P_4}\}}
|f_3(x)| dx\lesssim$$

$$\size_{3,\vec{T}_{3,4}}(f_3)
\sum_{\vec{P}\in\vec{T}}
|\langle f_1, \phi_{P_1}\rangle|
|\langle f_2, \phi_{P_2}\rangle|\lesssim
$$

$$\size_{3,\vec{T}_{3,4}}(f_3)
(\sum_{\vec{P}\in\vec{T}}
|\langle f_1, \phi_{P_1}\rangle|^2)^{1/2}
(\sum_{\vec{P}\in\vec{T}}
|\langle f_2, \phi_{P_2}\rangle|^2)^{1/2}\lesssim
$$

\begin{equation}
\size_{1, T_1}(f_1)\size_{1, T_2}(f_2)
\size_{3,\vec{T}_{3,4}}(f_3)|I_{\vec{T}}|
\end{equation}
and this ends the proof.

\end{document}